\documentclass{amsart}

\usepackage{amsmath}
\usepackage[latin1]{inputenc}
\usepackage{amsfonts}
\usepackage{amssymb}

\usepackage{amsthm}
\usepackage{amscd}

\usepackage[latin1]{inputenc}

\usepackage[T1]{fontenc}

\usepackage{amsmath}

\usepackage{amsfonts}

\usepackage{amssymb}

\usepackage{amsthm}

\usepackage{graphics}

\newcommand{\mk}{\medskip}

\newcommand{\ZZ}{\mathbb{Z}}

\newcommand{\CC}{\mathbb{C}}

\newcommand{\NN}{\mathbb{N}}

\newcommand{\yl}{15pt}

\newcommand{\ffbox}[1]{
\setbox9=\hbox{$\scriptstyle\overline{1}$}
\framebox[\yl][c]{\rule{0mm}{\ht9}${\scriptstyle #1}$}
}

\newcommand{\QQ}{\mathbb{Q}}

\newcommand{\Glie}{\mathfrak{g}}

\newcommand{\Yim}{\mathcal{Y}}

\newcommand{\Hlie}{\mathfrak{h}}

\newcommand{\demo}{\noindent {\it \small Proof:}\quad}

\newcommand{\U}{\mathcal{U}}

\newcommand{\Lo}{\mathcal{L}}

\newtheorem{thm}{Theorem}[section]

\newtheorem{defi}[thm]{Definition}

\newtheorem{cor}[thm]{Corollary}

\newtheorem{prop}[thm]{Proposition}

\newtheorem{lem}[thm]{Lemma}

\newtheorem{rem}[thm]{Remark}

\usepackage[all]{xy}
\title{On minimal affinizations of representations of quantum groups}

\author{David Hernandez}
\address{CNRS - UMR 8100 : Laboratoire de Math\'ematiques de Versailles, 45 avenue des Etats-Unis , Bat. Fermat, 78035 VERSAILLES, 
FRANCE}
\email{hernandez @ math . cnrs . fr}
\urladdr{http://www.math.uvsq.fr/\textasciitilde hernandez}

\begin{document}

\begin{abstract} In this paper we study minimal affinizations of representations of quantum groups (generalizations of Kirillov-Reshetikhin modules of quantum affine algebras introduced in \cite{Chari2}). We prove that all minimal affinizations in types $A$, $B$, $G$ are special in the sense of monomials. Although this property is not satisfied in general, we also prove an analog property for a large class of minimal affinization in types $C$, $D$, $F$. As an application, the Frenkel-Mukhin algorithm \cite{Fre2} works for these modules. For minimal affinizations of type $A$, $B$ we prove the thin property (the $l$-weight spaces are of dimension $1$) and a conjecture of \cite{nn1} (already known for type A). The proof of the special property is extended uniformly for more general quantum affinizations of quantum Kac-Moody algebras.

\vskip 4.5mm


\end{abstract}

\maketitle

\tableofcontents

\section{Introduction} 

In this paper $q\in\CC^*$ is fixed and is not a root of unity. 

Affine Kac-Moody algebras $\hat{\Glie}$ are infinite dimensional analogs of semi-simple Lie algebras $\Glie$, and have remarkable applications (see \cite{kac}). Their quantizations $\U_q(\hat{\Glie})$, called quantum affine algebras, have a very rich representation theory which has been intensively studied in mathematics and physics (see references in \cite{Cha2, dm} and in \cite{Cha, Fre, Naams, Nab} for various approaches). In particular Drinfeld \cite{Dri2} discovered that they can also be realized as quantum affinizations of usual quantum groups $\U_q(\Glie)\subset \U_q(\hat{\Glie})$. By using this new realization, Chari-Pressley \cite{Cha2} classified their finite dimensional representations.

Chari \cite{Chari2} introduced the notion of minimal affinizations of representations of quantum groups : starting from a simple representation $V$ of $\U_q(\hat{\Glie})$, an affinization of $V$ is a simple representation $\hat{V}$ of $\U_q(\hat{\Glie})$ such that $V$ is the head in the decomposition of $\hat{V}$ in simple $\U_q(\Glie)$-representations. Then one can define a partial ordering on the set of affinizations of $V$ and so a notion of minimal affinization for this ordering. For example the minimal affinizations of simple $\U_q(\Glie)$-modules of highest weight a multiple of a fundamental weight are the Kirillov-Reshetikhin modules which have been intensively studied in recent years (for example see \cite{kos, knh, kl, hkoty, k, c0, Nab, Nad, her06, cm3, fl} and references therein). An (almost) complete classification of minimal affinizations was done in \cite{Chari2, Cha7, Cha8, Cha9}.

The motivation to study minimal affinizations comes from physics : the affinizations of representations of quantum groups are important objects from the physical point of view, as stressed for example in \cite[Remark 4.2]{frere} and in the introduction of \cite{Chari2}. For example in the theory of lattice models in statistical mechanics, they are related to the problem of proving the integrability of the model : the point is to add spectral parameters to a solution of the related quantum Yang-Baxter equation (see \cite{Cha2}). A second example is related to the quantum particles of the affine Toda field theory (see \cite{bele, pdor}) which correspond to simple finite dimensional representations of quantum affine algebras.

In the present paper we prove new results on the structure of minimal affinizations, in particular in the light of recent developments in the representation theory of quantum affine algebras. 

A particular class of finite dimensional representations, called special modules \cite{Nab}, attracted much attention as Frenkel-Mukhin \cite{Fre2} proposed an algorithm which gives their $q$-character (analog of the usual character adapted to the Drinfeld realization and introduced by Frenkel-Reshetikhin \cite{Fre} : they encode a certain decomposition of representations in so called $l$-weight spaces or pseudo weight spaces). For example the Kirillov-Reshetikhin modules \cite{Nab, Nad, her06} are special (this is the crucial point of the proof of the Kirillov-Reshetikhin conjecture). A dual class of modules called antispecial modules is introduced in the present paper (antispecial does not mean the opposite of special), and an analog of the Frenkel-Mukhin algorithm gives their $q$-character.

In the present paper we prove that minimal affinizations in type $A$, $B$, $G$ are special and antispecial. We get counter examples for other types, but we prove in type $C$, $D$, $F$ that a large class of minimal affinizations are special or antispecial. In particular the Frenkel-Mukhin algorithm works for these modules. As an application, we prove that minimal affinizations of type $A$ and $B$ are thin (the $l$-weight spaces are of dimension $1$). We also get the special property for analog simple modules of quantum affinizations of some non necessarily finite quantum Kac-Moody algebras.

 In the proofs of the present paper, the crucial steps include technics developed in \cite{her06} to prove the Kirillov-Reshetikhin conjecture and in \cite{small} to solve the Nakajima's smallness problem. The general idea is to prove simultaneously the special property and the thin property by induction on the highest weight of the minimal affinizations. This allows to use the elimination theorem \cite{her06} which leads to eliminate some monomials in the $q$-character of simple modules.

 Nakajima first conjectured the existence of such large classes of special modules for simply-laced cases (see \cite{Nab}), and the existence of a large class of special minimal affinizations was conjectured by Mukhin in a conversation with the author in the conference "Representations of Kac-Moody Algebras and Combinatorics" at Banff in March 2005. 

In some situations, the properties are already known or can be proved directly from already known explicit formulas. Indeed, for Kirillov-Reshetikhin modules the special property was proved in \cite{Nad} (simply-laced case) and in \cite{her06} (non simply-laced case). So for Kirillov-Reshetikhin modules in classical types, the explicit formulas in \cite{kos, knh} are satisfied (the formulas for fundamental representations are given in \cite{ks}) and we can get the properties directly from them. General formulas and the thin properties were proved for irreducible
tame modules, which include minimal affinizations, for Yangians of
type $A$ \cite{che1, che, nt}. (The author was told by Nakajima that the
same results hold for quantum affine algebras of type $A$ by \cite{v}.) See
also \cite{Fre3} for the cases of minimal affinizations, which are
evaluation representations in type $A$.

 

Explicit formulas are also available for twisted yangians in classical types \cite{mo, naz}. But the author did not find in the literature a proof of the correspondence between quantum affine algebras and twisted (or non simply laced) yangians. 

In general no explicit formulas for $q$-characters of quantum affine algebras are available, so our proofs use direct arguments without explicit formulas and are independent of previous results on yangians. In particular this allows to extend uniformly our arguments to previously unknown situations (like type $B,C,D,G_2,F_4$), and to more general quantum affinizations which are not necessarily quantum affine algebras.


For quantum affine algebras in classical types, explicit conjectural formulas \cite{nn1, nn2, nn3} are available for a large class of representations including many minimal affinizations (all of them for type $A$; see \cite{kos, knh} for more general formulas). In types $A$, $B$, the results proved in the present paper imply \cite[Conjecture 2.2]{nn1} for these minimal affinizations. The author did not find in the literature a proof of this result, except for type $A$ as explained above. The main subject of the present paper is minimal affinizations and so we give a proof of \cite[Conjecture 2.2]{nn1} in this case. But it is possible to prove \cite[Conjecture 2.2]{nn1} for more general representations by using a variation of this proof (this and \cite[Conjecture 2.2]{nn1} in types $C$, $D$ will be discussed in a separate publication).

The results of \cite{nt, ks, kos} and of \cite[Conjecture 2.2]{nn1} (and thin property as their consequence) were explained to the author by Nakajima in an early stage of this research, June 2005.



Let us describe the organization of the present paper. In section \ref{un} we give some backgrounds on the representation theory of quantum affine algebras. 
In section \ref{deux} we recall the definition of minimal affinizations and state the main results of the paper.
In section \ref{trois} we give preliminary results, including results from \cite{small} and discussion about an involution of $\U_q(\hat{\Glie})$. In section \ref{proofmainres} we prove the main result of the paper. 
In section \ref{lastsec} we explain the proof of \cite[Conjecture 2.2]{nn1} for minimal affinizations in types $A$, $B$, we state additional results (Theorem \ref{moregen}) for more general quantum affinizations, and we discuss possible further developments, in particular on generalized induction systems involving minimal affinizations.

{\bf Acknowledgments :} The author is very grateful Evgeny Mukhin for encouraging him to study minimal affinizations in the continuation of the proof of the Kirillov-Reshetikhin conjecture, to Hiraku Nakajima for useful comments and references in an early stage of this research, and to Vyjayanthi Chari, Maxim Nazarov, Alexander Molev, Michela Varagnolo for useful comments and references. A part of this paper was written as the author gave lectures in the East China Normal University in Shanghai, he would like to thank Naihong Hu for the invitation.

\section{Background}\label{un}

\subsection{Cartan matrix and quantized Cartan matrix} Let $C=(C_{i,j})_{1\leq i,j\leq n}$ be a Cartan matrix of 
finite type. We denote 
$I=\{1,\cdots,n\}$. $C$ is symmetrizable : there is a matrix $D=\text{diag}(r_1,\cdots,r_n)$ ($r_i\in\NN^*$)\label{ri} such 
that $B=DC$\label{symcar} is symmetric. In particular if $C$ is symmetric then $D=I_n$ (simply-laced case). 

\noindent We consider a realization $(\Hlie, \Pi, \Pi^{\vee})$ of $C$ (see \cite{bou, kac}): $\Hlie$ is a $n$ dimensional $\QQ$-vector space, $\Pi=\{\alpha_1,\cdots ,\alpha_n\}\subset \Hlie^*$ (set of the simple roots) and $\Pi^{\vee}=\{\alpha_1^{\vee},\cdots ,\alpha_n^{\vee}\}\subset \Hlie$ (set of simple coroots) are set such that for $1\leq i,j\leq n$, $\alpha_j(\alpha_i^{\vee})=C_{i,j}$.
Let $\Lambda_1,\cdots ,\Lambda_n\in\Hlie^*$ (resp. $\Lambda_1^{\vee},\cdots ,\Lambda_n^{\vee}\in\Hlie$) be the the fundamental weights (resp. coweights) : $\Lambda_i(\alpha_j^{\vee})=\alpha_i(\Lambda_j^{\vee})=\delta_{i,j}$ where $\delta_{i,j}$ is $1$ if $i=j$ and $0$ otherwise. Denote $P=\{\lambda \in\Hlie^*|\forall i\in I, \lambda(\alpha_i^{\vee})\in\ZZ\}$ the set of weights and $P^+=\{\lambda \in P|\forall i\in I, \lambda(\alpha_i^{\vee})\geq 0\}$ the set of dominant weights. For example we have $\alpha_1,\cdots ,\alpha_n\in P$ and $\Lambda_1,\cdots ,\Lambda_n\in P^+$. Denote $Q={\bigoplus}_{i\in I}\ZZ \alpha_i\subset P$ the root lattice and $Q^+={\sum}_{i\in I}\NN \alpha_i\subset Q$. For $\lambda,\mu\in \Hlie^*$, denote $\lambda \geq \mu$ if $\lambda-\mu\in Q^+$. Let $\nu:\Hlie^*\rightarrow \Hlie$ linear such that for all $i\in I$ we have $\nu(\alpha_i)=r_i\alpha_i^{\vee}$. For $\lambda,\mu\in\Hlie^*$, $\lambda(\nu(\mu))=\mu(\nu(\lambda))$. We use the enumeration of vertices of \cite{kac}.

\noindent We denote $q_i=q^{r_i}$ and for $l\in\ZZ, r\geq 0, m\geq m'\geq 0$ we define in $\ZZ[q^{\pm}]$ :
$$[l]_q=\frac{q^l-q^{-l}}{q-q^{-1}}\text{ , }[r]_q!=[r]_q[r-1]_q\cdots [1]_q\text{ ,
}\begin{bmatrix}m\\m'\end{bmatrix}_q=\frac{[m]_q!}{[m-m']_q![m']_q!}.$$

\noindent For $a,b\in\ZZ$, we denote $q^{a+b\ZZ}=\{q^{a+br}|r\in\ZZ\}$ and $q^{a+b\NN}=\{q^{a+br}|r\in\ZZ,r\geq 0\}$.

\noindent Let $C(z)$ be the quantized Cartan matrix defined by ($i\neq j\in I$): 
$$C_{i,i}(z)=z_i+z_i^{-1}\text{ ,
}C_{i,j}(z)=[C_{i,j}]_z.$$ 
$C(z)$ is invertible (see \cite{Fre}). We denote by $\tilde{C}(z)$ the inverse matrix of
$C(z)$ and by $D(z)$ the diagonal matrix such that for $i,j\in I$, $D_{i,j}(z)=\delta_{i,j}[r_i]_z$.

\subsection{Quantum algebras}\label{qkma}

\subsubsection{Quantum groups}

\begin{defi} The quantum group $\U_q(\Glie)$ is the $\CC$-algebra with generators $k_i^{\pm 1}$, $x_i^{\pm}$ ($i\in I$) and 
relations: 
$$k_ik_j=k_jk_i\text{ , } k_ix_j^{\pm}=q_i^{\pm C_{i,j}}x_j^{\pm}k_i,$$
$$[x_i^+,x_j^-]=\delta_{i,j}\frac{k_i-k_i^{-1}}{q_i-q_i^{-1}},$$
$${\sum}_{r=0\cdots  1-C_{i,j}}(-1)^r\begin{bmatrix}1-C_{i,j}\\r\end{bmatrix}_{q_i}(x_i^{\pm})^{1-C_{i,j}-r}x_j^{\pm}(x_i^{\pm})^r=0 \text{ (for $i\neq j$)}.$$
\end{defi}

\noindent This algebra was introduced independently by Drinfeld \cite{Dri1} and Jimbo \cite{jim}. It is remarkable that 
one can define a Hopf algebra structure on $\U_q(\Glie)$ by : 
$$\Delta(k_i)=k_i\otimes k_i,$$ 
$$\Delta(x_i^+)=x_i^+\otimes 1 + k_i\otimes x_i^+\text{ , }\Delta(x_i^-)=x_i^-\otimes 
k_i^{-1}+ 1\otimes x_i^-,$$ 
$$S(k_i)=k_i^{-1}\text{ , }S(x_i^+)=-x_i^+k_i^{-1}\text{ , }S(x_i^-)=-k_ix_i^-,$$ 
$$\epsilon(k_i)=1\text{ , }\epsilon(x_i^+)=\epsilon(x_i^-)=0.$$

\noindent Let $\U_q(\Hlie)$ be the commutative subalgebra of $\U_q(\Glie)$ generated by the $k_i^{\pm 1}$ ($i\in I$).

\noindent For $V$ a $\U_q(\Hlie)$-module and $\omega\in P$ we denote by $V_{\omega}$ the weight space of weight 
$\omega$ : 
$$V_{\omega}=\{v\in V|\forall i\in I, k_i.v=q_i^{\omega(\alpha_i^{\vee})}v\}.$$ 
In particular we have $x_i^{\pm}.V_{\omega}\subset V_{\omega \pm \alpha_i}$.

\noindent We say that $V$ is $\U_q(\Hlie)$-diagonalizable if $V={\bigoplus}_{\omega\in P} V_{\omega}$ (in particular $V$ is of type $1$).

\noindent For $V$ a finite dimensional $\U_q(\Hlie)$-diagonalizable module we define the usual character 
$$\chi(V)={\sum}_{\omega\in P}\text{dim}(V_{\omega})e(\omega)\in\mathcal{E}=\bigoplus_{\omega\in P}\ZZ.e(\omega).$$ 

\subsubsection{Quantum loop algebras} We will use the second realization (Drinfeld realization) of the quantum loop 
algebra $\U_q(\Lo\Glie)$ (subquotient of the quantum affine algebra $\U_q(\hat{\Glie})$) :

\begin{defi}\label{defiaffi} $\U_q(\Lo\Glie)$ is the algebra with
generators $x_{i,r}^{\pm}$ ($i\in I, r\in\ZZ$), $k_i^{\pm 1}$ ($i\in I$), $h_{i,m}$ ($i\in I, m\in\ZZ-\{0\}$) and the
following relations ($i,j\in I, r,r'\in\ZZ, m,m' \in\ZZ-\{0\}$): 
$$\label{afcart}[k_i,k_j]=[k_i,h_{j,m}]=[h_{i,m},h_{j,m'}]=0,$$
$$k_ix_{j,r}^{\pm}=q_i^{\pm C_{i,j}}x_{j,r}^{\pm}k_i,$$
$$[h_{i,m},x_{j,r}^{\pm}]=\pm \frac{1}{m}[mB_{i,j}]_qx_{j,m+r}^{\pm},$$
$$[x_{i,r}^+,x_{j,r'}^-]=\delta_{i,j}\frac{\phi^+_{i,r+r'}-\phi^-_{i,r+r'}}{q_i-q_i^{-1}},$$
$$x_{i,r+1}^{\pm}x_{j,r'}^{\pm}-q^{\pm B_{i,j}}x_{j,r'}^{\pm}x_{i,r+1}^{\pm}=q^{\pm B_{i,j}}x_{i,r}^{\pm}x_{j,r'+1}^{\pm}-x_{j,r'+1}^{\pm}x_{i,r}^{\pm},$$
$${\sum}_{\pi\in\Sigma_s}{\sum}_{k=0\cdots s}(-1)^k\begin{bmatrix}s\\k\end{bmatrix}_{q_i}x_{i,r_{\pi(1)}}^{\pm}\cdots x_{i,r_{\pi(k)}}^{\pm}x_{j,r'}^{\pm}x_{i,r_{\pi(k+1)}}^{\pm}\cdots x_{i,r_{\pi(s)}}^{\pm}=0,$$
where the last relation holds for all $i\neq j$, $s=1-C_{i,j}$, all sequences of integers $r_1,\cdots ,r_s$. $\Sigma_s$ is
the symmetric group on $s$ letters. For $i\in I$ and $m\in\ZZ$, $\phi_{i,m}^{\pm}\in \U_q(\Lo\Glie)$ is determined
by the formal power series in $\U_q(\Lo\Glie)[[z]]$ (resp. in $\U_q(\Lo\Glie)[[z^{-1}]]$): 
$${\sum}_{m\geq 0}\phi_{i,\pm m}^{\pm}z^{\pm m}=k_i^{\pm}\text{exp}(\pm(q-q^{-1}){\sum}_{m'\geq 1}h_{i,\pm m'}z^{\pm m'}),$$ 
and $\phi_{i,\mp m}^{\pm}=0$ for $m>0$. \end{defi}

\noindent $\U_q(\Lo\Glie)$ has a Hopf algebra structure (from the Hopf algebra structure of $\U_q(\hat{\Glie})$).

\noindent For $J\subset I$ we denote by $\U_q(\Lo\Glie_J)\subset\U_q(\Lo\Glie)$ the subalgebra generated by the $x_{i,m}^{\pm}$, $h_{i,m}$, $k_i^{\pm 1}$ for $i\in J$. $\U_q(\Lo\Glie_J)$ is a quantum loop algebra associated to the semi-simple Lie algebra $\Glie_J$ of Cartan matrix $(C_{i,j})_{i,j\in J}$. For example for $i\in I$, we denote $\U_q(\Lo\Glie_i) = \U_q(\Lo\Glie_{\{i\}})\simeq \U_{q_i}(\Lo sl_2)$.

\noindent The subalgebra of $\U_q(\Lo\Glie)$ generated by the $h_{i,m}, k_i^{\pm 1}$ (resp. by the $x_{i,r}^{\pm}$) is denoted by $\U_q(\Lo\Hlie)$ (resp. $\U_q(\Lo\Glie)^{\pm}$).

\subsection{Finite dimensional representations of quantum loop algebras} Denote by 
$\text{Rep}(\U_q(\Lo\Glie))$ the Grothendieck ring of (type $1$) finite 
dimensional representations of $\U_q(\Lo\Glie)$.

\subsubsection{Monomials and $q$-characters}\label{defimono} Let $V$ be a representation in $\text{Rep}(\U_q(\Lo\Glie))$. The subalgebra $\U_q(\Lo\Hlie)\subset\U_q(\Lo\Glie)$ is commutative, so we have :
$$V={\bigoplus}_{\gamma=(\gamma_{i,\pm m}^{\pm})_{i\in I, m\geq 0}}V_{\gamma},$$
$$\text{where : }V_{\gamma}=\{v\in V|\exists p\geq 0, \forall i\in I, m\geq 0, (\phi_{i,\pm m}^{\pm}-\gamma_{i,\pm m}^{\pm})^p.v=0\}.$$
The $\gamma=(\gamma_{i,\pm m}^{\pm})_{i\in I, m\geq 0}$ are called $l$-weights (or pseudo-weights) and the $V_{\gamma}\neq \{0\}$ are called $l$-weight spaces (or pseudo-weight spaces) of $V$. One can prove \cite{Fre} that $\gamma$ is necessarily of the form :
\begin{equation}\label{formfin}{\sum}_{m\geq 0}\gamma_{i,\pm m}^{\pm}u^{\pm 
m}=q_i^{\text{deg}(Q_i)-\text{deg}(R_i)}\frac{Q_i(uq_i^{-1})R_i(uq_i)}{Q_i(uq_i)R_i(uq_i^{-1})},
\end{equation}
where $Q_i,R_i\in \CC(u)$ satisfy $Q_i(0)=R_i(0)=1$. The Frenkel-Reshetikhin $q$-characters morphism $\chi_q$ \cite{Fre} encodes the $l$-weights $\gamma$ (see also \cite{kn}). It is an injective ring morphism : 
$$\chi_q:\text{Rep}(\U_q(\Lo\Glie))\rightarrow \ZZ[Y_{i,a}^{\pm}]_{i\in I, a\in\CC^*}$$
defined by
$$\chi_q(V) = {\sum}_{\gamma}\text{dim}(V_{\gamma})m_{\gamma},$$
where :
$$m_{\gamma}={\prod}_{i\in I, a\in\CC^*}Y_{i,a}^{q_{i,a}-r_{i,a}},$$ 
$$Q_i(u)={\prod}_{a\in\CC^*}(1-ua)^{q_{i,a}}\text{ , }R_i(u)={\prod}_{a\in\CC^*}(1-ua)^{r_{i,a}}.$$
The $m_{\gamma}$ are called monomials (they are analogs of weight). We denote by $A$ the set of monomials of $\ZZ[Y_{i,a}^{\pm}]_{i\in I, a\in\CC^*}$. For an $l$-weight $\gamma$, we denote $V_{\gamma}=V_{m_{\gamma}}$. We will also use the notation $i_r^p=Y_{i,q^r}^p$ for $i\in I$ and $r,p\in\ZZ$.

\noindent For $J\subset I$, $\chi_q^J$ is the morphism of $q$-characters for $\U_q(\Lo\Glie_J)\subset\U_q(\Lo\Glie)$. 
\\For a $m$ monomial we denote $u_{i,a}(m)\in\ZZ$ such that $m = {\prod}_{i\in I, a\in\CC^*}Y_{i,a}^{u_{i,a}(m)}$. We also denote $\omega(m) = {\sum}_{i\in I,a\in\CC^*} u_{i,a}(m)\Lambda_i$, $u_i(m)=\sum_{a\in\CC^*}u_{i,a}(m)$ and $u(m)=\sum_{i\in I}u_i(m)$. $m$ is said to be 
$J$-dominant if for all $j\in J, a\in\CC^*$ we have $u_{j,a}(m)\geq 0$. An $I$-dominant monomials is said to be dominant.

\noindent Observe that $\chi_q, \chi_q^J$ can also be defined for finite dimensional $\U_q(\Lo \Hlie)$-modules in the same way.

\noindent In the following for $V$ a finite dimensional $\U_q(\Lo\Hlie)$-module, we denote by $\mathcal{M}(V)$ the set 
of monomials occurring in $\chi_q(V)$. 

\noindent For $i\in I, a\in\CC^*$ we set : 
\begin{equation}\label{aia}
\begin{split}
A_{i,a}=&Y_{i,aq_i^{-1}}Y_{i,aq_i}\prod_{\{j|C_{j,i}=-1\}}Y_{j,a}^{-1}
\\&\times\prod_{\{j|C_{j,i}=-2\}}Y_{j,aq^{-1}}^{-1}Y_{j,aq}^{-1}\prod_{\{j|C_{j,i}=-3\}}Y_{j,aq^2}^{-1}Y_{j,a}^{-1}Y_{j,aq^{-2}}^{-1}.
\end{split}
\end{equation}

\noindent As the $A_{i,a}^{-1}$ are algebraically independent \cite{Fre} (because $C(z)$ is invertible), for $M$ a product of $A_{i,a}^{-1}$ we can define $v_{i,a}(M)\geq 0$ by $M={\prod}_{i\in I, a\in\CC^*} A_{i,a}^{-v_{i,a}(m)}$. We put $v_i(M)={\sum}_{a\in\CC^*}v_{i,a}(M)$ and $v(M)={\sum}_{i\in I}v_i(M)$.

\noindent For $\lambda\in - Q^+$ we set $v(\lambda)=-\lambda(\Lambda_1^{\vee}+\cdots +\Lambda_n^{\vee})$. For $M$ a 
product of $A_{i,a}^{-1}$, we have $v(M)=v(\omega(\lambda))$.

\noindent For $m,m'$ two monomials, we write $m'\leq m$ if $m'm^{-1}$ is product of $A_{i,a}^{-1}$. 

\begin{defi}\label{monomrn}\cite{Fre2} A monomial $m\in A-\{1\}$ is said to be right-negative if for all $a\in\CC^*$, 
for $L=\text{max}\{l\in\ZZ|\exists i\in I, u_{i,aq^L}(m)\neq 0\}$ we have $\forall j\in I$, $u_{j,aq^L}(m)\leq 0$.\end{defi}

\noindent Observe that a right-negative monomial is not dominant. We can also define left-negative monomials by replacing $\text{max}$ by $\text{min}$ in the formula of $L$ in Definition \ref{monomrn}.

\begin{lem}\cite{Fre2}\label{rn} 1) For $i\in I, a\in\CC^*$, $A_{i,a}^{-1}$ is right-negative.

2) A product of right-negative monomials is right-negative.

3) If $m$ is right-negative, then $m'\leq m$ implies that $m'$ is right-negative.\end{lem}

We have the same results by replacing right-negative by left-negative.

\noindent For $J\subset I$ and $Z\in \Yim$, we denote $Z^{\rightarrow J}$ the element of $\Yim$ obtained from $Z$ by putting $Y_{j,a}^{\pm 1}=1$ for $j\notin J$.

\noindent Let $\beta : \ZZ[Y_{j,b}^{\pm}]_{j\in I,b\in\CC^*}\rightarrow \mathcal{E}$ be the ring morphism such that $\beta(m)=e(\omega(m))$. 

\begin{prop}\label{restrcar}\cite[Theorem 3]{Fre} For $V\in\text{Rep}(\U_q(\Lo\Glie))$, let $\text{Res}(V)$ be the restricted 
$\U_q(\Glie)$-module. We have $(\beta\circ\chi_q)(V)=\chi(\text{Res}(V))$.\end{prop}

\subsubsection{$l$-highest weight representations}

\noindent The irreducible finite dimensional $\U_q(\Lo\Glie)$-modules have been classified by Chari-Pressley. They are parameterized by dominant monomials : 

\begin{defi}\label{lhigh} A $\U_q(\Lo\Glie)$-module $V$ is said to be of $l$-highest weight $m\in A$ if there is $v\in V_m$ such that $V=\U_q(\Lo\Glie)^-.v$ and $\forall i\in I, r\in\ZZ, x_{i,r}^+.v=0$.\end{defi}

\noindent For $m\in A$, there is a unique simple module $L(m)$ of $l$-highest weight $m$. 

\begin{thm}\cite[Theorem 12.2.6]{Cha2} The dimension of $L(m)$ is finite if and only if $m$ is dominant.\end{thm}

For $i\in I$, $a\in\CC^*$, $k\geq 0$ we denote $X_{k,a}^{(i)}={\prod}_{k'\in\{1,\cdots ,k\}}Y_{i,aq_i^{k-2k'+1}}$. 

\begin{defi} The Kirillov-Reshetikhin modules are the $W_{k,a}^{(i)} = L(X_{k,a}^{(i)})$.\end{defi}

We denote by $W_{0,a}^{(i)}$ the trivial representation (it is of dimension $1$). For $i\in I$ and $a\in\CC^*$, $W_{1,a}^{(i)}$ is called a fundamental representation and is denoted by $V_i(a)$ (in the case $\Glie = sl_2$ we simply write $W_{k,a}$ and $V(a)$). 

\noindent For $\Glie = sl_2$, the monomials $m_1 = X_{k_1,a_1}$, $m_2 = X_{k_2,a_2}$ are said to be in special position if the monomial $m_3={\prod}_{a\in\CC^*}Y_a^{\text{max}(u_a(m_1),u_a(m_2))}$ is of the form $m_3 = X_{k_3,a_3}$ and $m_3\neq m_1$, $m_3\neq m_2$. A normal writing of a dominant monomial $m$ is a product decomposition $m={\prod}_{i=1,\cdots ,L}X_{k_l,a_l}$ such that for $l\neq l'$, $X_{k_l,a_l}$, $X_{k_{l'},a_{l'}}$ are not in special position. Any dominant monomial has a unique normal writing up to permuting the monomials (see \cite[Section 12.2]{Cha2}).

\noindent It follows from the study of the representations of $\U_q(\Lo sl_2)$ in \cite{Cha0, Cha, Fre} that :

\begin{prop}\label{aidesldeux} Suppose that $\Glie=sl_2$.

(1) $W_{k,a}$ is of dimension $k+1$ and :
$$\chi_q(W_{k,a})=X_{k,a}(1+A_{aq^{k}}^{-1}(1+A_{aq^{k-2}}^{-1}(1+\cdots (1+A_{aq^{2-k}}^{-1}))\cdots ).$$

(2) $V(aq^{1-k})\otimes V(aq^{3-k})\otimes \cdots \otimes V(aq^{k-1})$ is of $q$-character :
$$X_{k,a}(1+A_{aq^k}^{-1})(1+A_{aq^{k-2}}^{-1})\cdots (1+A_{aq^{2-k}}^{-1}).$$
In particular all $l$-weight spaces of the tensor product are of dimension $1$.

(3) For $m$ a dominant monomial and $m = X_{k_1,a_1}\cdots X_{k_l,a_l}$ a normal writing we have :
$$L(m)\simeq W_{k_1,a_1}\otimes \cdots \otimes W_{k_l,a_l}.$$
\end{prop}

\subsubsection{Special modules and complementary reminders}\label{compred}

\begin{defi} For $m\in A$ let $D(m)$ be the set of monomials $m'\in A$ such
that there are $m_0=m, m_1, \cdots, m_N=m'\in A$ satisfying for all $j\in\{1,\cdots ,N\}$ : 
\begin{enumerate}
\item $m_j=m_{j-1}A_{i_j,a_1q_{i_j}}^{-1}\cdots A_{i_j,a_{r_j}q_{i_j}}^{-1}$ where $i_j\in I$, $r_j\geq 1$ and $a_1,\cdots, a_{r_j}\in\CC^*$,

\item for $1\leq r\leq r_j$, $u_{i_j,a_r}(m_{j-1})\geq
|\{r'\in\{1,\cdots ,r_j\}|a_{r'}=a_r\}|$ where $r_j,i_j,a_r$ are as in condition (1).
\end{enumerate}\end{defi}

\noindent For all $m'\in D(m)$, $m'\leq m$. Moreover if $m'\in D(m)$, then $(D(m')\subset D(m))$.

\begin{thm}\label{inducb}\cite[Theorem 5.21]{her07} For $V\in\text{Mod}(\U_q(\hat{\Glie}))$ be a $l$-highest weight module of highest monomial $m$, we have $\mathcal{M}(V)\subset D(m)$.\end{thm}

\noindent In particular for all $m'\in\mathcal{M}(V)$, we have $m'\leq m$ and the $v_{i,a}(m'm^{-1}), v(m'm^{-1})\geq 0$ are well-defined. As a direct consequence of Theorem \ref{inducb}, we also have : 

\begin{lem}\label{fundless} For $i\in I, a\in\CC^*$, we have $(\chi_q(V_i(a))-Y_{i,a})\in\ZZ[Y_{j,aq^l}^{\pm}]_{j\in I, l > 0}$.\end{lem}

\noindent This result was first proved in \cite[Lemma 6.1, Remark 6.2]{Fre2}.

\noindent A monomial $m$ is said to be antidominant if for all $i\in I, a\in\CC^*$, $u_{i,a}(m)\leq 0$. 

\begin{defi} A $\U_q(\Lo\Glie)$-module is said to be special (resp. antispecial) if his $q$-character has a unique dominant (resp. antidominant) monomial.\end{defi}

\noindent The notion of special module was introduced in \cite{Nab}. It is of particular importance because an algorithm of Frenkel-Mukhin \cite{Fre2} gives the $q$-character of special modules. It is easy to write a similar algorithm for antispecial modules from the Frenkel-Mukhin algorithm (for example it suffices to use the involution studied in section \ref{secinvol}).

\noindent Observe that a special (resp. antispecial) module is a simple $l$-highest weight module. But in general all simple $l$-highest weight module are not special. The following result was proved in \cite{Nab, Nad} for simply-laced types, and in full generality in \cite{her06} (see \cite{Fre2} for previous results).

\begin{thm}\cite[Theorem 4.1, Lemma 4.4]{her06}\label{formerkr} The Kirillov-Reshetikhin modules are special. Moreover for $m\in \mathcal{M}(W_{k,a}^{(i)})$, $m\neq X_{k,a}^{(i)}$ implies $m\leq X_{k,a}^{(i)}A_{i,aq_i^k}^{-1}$.\end{thm}

 Define 
$$\mu_J^I:\ZZ[(A_{j,a}^{\pm})^{\rightarrow (J)}]_{j\in J,a\in\CC^*}\rightarrow \ZZ[A_{j,a}^{\pm}]_{j\in J,a\in\CC^*}$$ 
the ring morphism such that $\mu_J^I((A_{j,a}^{\pm})^{\rightarrow (J)})=A_{j,a}^{\pm}$. For $m$ $J$-dominant, denote by $L^J(m^{\rightarrow (J)})$ the simple $\U_q(\Lo\Glie_J)$-module of $l$-highest weight $m^{\rightarrow (J)}$. Define :
$$L_J(m)=m \mu_J^I((m^{\rightarrow (J)})^{-1}\chi_q^J(L^J(m^{\rightarrow (J)}))).$$
We have :

\begin{prop}\label{jdecomp}\cite{small} For a representation $V\in\text{Rep}(\U_q(\Lo \Glie))$ and $J\subset I$, there is unique
decomposition in a finite sum : 
\begin{equation}\label{formdecomp}\chi_q(V)=\sum_{m'\text{ $J$-dominant}}\lambda_J(m')L_J(m').\end{equation}
Moreover for all $m'$ $J$-dominant we have $\lambda_J(m')\geq 0$.\end{prop}

\begin{rem}\label{process} Let $m$ be a dominant monomial and $m'\in\mathcal{M}(L(m))$ a $J$-dominant monomial such that there are no $m'' > m'$ satisfying $m''\in\mathcal{M}(m)$ and $m'$ appears in $L_J(m'')$. Then from Proposition \ref{jdecomp} the monomials of $L_J(m')$ are in $\mathcal{M}(L(m))$. It gives inductively from $m$ a set of monomial occurring in $\chi_q(L(m))$.\end{rem}

\subsubsection{Thin modules}

\begin{defi}\cite{small} A $\U_q(\Lo\Glie)$-module $V$ is said to be thin if his $l$-weight spaces are of dimension $1$.\end{defi}

For example for $\Glie$ of type $A$, $B$, $C$, all fundamental representations are thin (it can be established directly from the formulas in \cite{ks}; this thin property was observed and proved with a different method in \cite[Theorem 3.5]{her05}; see also \cite{cm2}). 

\noindent Observe that it follows from \cite[Section 8.4]{her02} that for $\Glie$ of type $G_2$, all fundamental representations are thin. For $\Glie$ of type $F_4$, the fundamental representations corresponding to $i=1$ and $i=4$ are thin, but the fundamental representations corresponding to $i=2$ or $i=3$ are not thin (see \cite{her05}).

\noindent For type $D$, it is known that fundamental representations are not necessarily thin : for example for $\Glie$ of type $D_4$, the fundamental representations $V_2(a)$ has an $l$-weight space of dimension $2$. Explicit formulas for the $q$-character of fundamental representation in type $D$ are given in \cite{ks} (the thin fundamental representations of type $D$ are also characterized in \cite{cm2}; see also remark \ref{remnak} bellow for a general statement).

For $m\in\ZZ[Y_{i,a}]_{i\in I, a\in\CC^*}$ a dominant monomial, the standard module $M(m)$ is defined as the tensor product : 
$$M(m)=\bigotimes_{a\in(\CC^*/q^{\ZZ})}(\cdots\otimes(\bigotimes_{i\in I}V_i(aq)^{\otimes u_{i,aq}(m)})\otimes(\bigotimes_{i\in I}V_i(aq^2)^{\otimes u_{i,aq^2}(m)})\otimes\cdots).$$
It is well-defined as for $i,j\in I$ and $a\in\CC^*$ we have $V_i(a)\otimes V_j(a)\simeq V_j(a)\otimes V_i(a)$ and for $a'\notin aq^{\ZZ}$, we have $V_i(a)\otimes V_j(a')\simeq V_j(a')\otimes V_i(a)$. Observe that fundamental representations are particular cases of standard modules.

As a direct corollary of a result of Nakajima, there is the following result for simply-laced types :

\begin{cor} We suppose that $\Glie$ is simply-laced. Let $m=\prod_{i\in I}Y_{i,aq^{\phi_i}}^{w_i}$ where $a\in\CC^*$, $w_i\geq 0$ and $\phi_i\in\ZZ$ satisfies $(C_{i,j}<0\Rightarrow |\phi_i - \phi_j|=1)$. Then the standard module $M(m)$ is thin if and only if it is simple as a $\U_q(\Glie)$-module.\end{cor}

\demo It follows from \cite[Proposition 3.4]{Naexa} that in this situation the number of monomials in $\chi_q(M(m)))$ is equal to the dimension of the simple $\U_q(\Glie)$-module of highest weight $\sum_{i\in I}w_i\Lambda_i$.\qed

Observe that this result is false for not simply-laced $\Glie$ (for example there is a thin fundamental representation for type $G_2$ which is not simple as a $\U_q(\Glie)$-module, see \cite[Section 8.4]{her02}).

The following remark was communicated to the author by Nakajima :

\begin{rem}\label{remnak} In particular for $\Glie$ simply-laced, a fundamental representation is thin if and only if the corresponding coefficient of the highest root is $1$ (this point is also a trivial consequence of previously known results, for example the geometric construction \cite{Naams}).\end{rem}

We got also the following example :

\begin{prop}\cite[Proposition 6.6]{small} Let $\Glie=sl_{n+1}$ and consider a monomial of the form $m = Y_{i_1,aq^{l_1}}Y_{i_2,aq^{l_2}}\cdots Y_{i_R,aq^{l_R}}$ where $R\geq 0$, $i_1,i_2,\cdots, i_R\in I$, $l_1,l_2,\cdots , l_R\in \ZZ$ satisfying for all $1\leq r\leq R-1$, $l_{r+1}-l_r \geq i_r+i_{r+1}$. Then $L(m)$ is thin.
\end{prop}

\section{Minimal affinizations and main results}\label{deux} In this section we recall the definition of minimal affinizations and their classification in regular cases. Then we state the main results which are proved in other sections. 
\subsection{Definitions \cite{Chari2}}

\begin{defi} For $V$ a simple finite dimensional $\U_q(\Glie)$-module, a simple finite dimensional $\U_q(\Lo\Glie)$-module $L(m)$ is said to be an affinization of $V$ if $\omega(m)$ is the highest weight of $V$.\end{defi}

\noindent For $V$ a $\U_q(\Glie)$-module and $\lambda\in P^+$, denote by $m_{\lambda}(V)$ the multiplicity in $V$ of the simple $\U_q(\Glie)$-module of highest weight $\lambda$.

\noindent Two affinizations are said to equivalent if they are isomorphic as $\U_q(\Glie)$-modules. Denote by $\mathcal{Q}_{V}$ the equivalence classes of affinizations of $V$ and for $L$ an affinization of $V$ denote by $[L]\in\mathcal{Q}_V$ its classes. For $[L], [L']\in \mathcal{Q}_V$, we write $[L]\leq [L']$ if and only if for all $\mu\in P^+$, either

(i) $m_{\mu}(L)\leq m_{\mu}(L')$,

(ii) $\exists \nu>\mu$ such that $m_{\nu}(L)< m_{\mu}(L')$.

\begin{prop} $\leq$ defines a partial ordering on $\mathcal{Q}_V$.\end{prop}

\begin{defi} A minimal affinizations of $V$ is a minimal element of $\mathcal{Q}_V$ for the partial ordering.\end{defi}

\begin{rem} For $\Glie = sl_{n+1}$, we have evaluation morphisms $\U_q(\Lo \Glie)\rightarrow \U_q(\Glie)$ denoted by $ev_a$ and $ev^a$ (for $a\in \CC^*$) and in particular a minimal affinization $L$ of $V$ is isomorphic to $V$ as a $\U_q(\Glie)$-module.\end{rem}

\subsection{Classification}

The minimal affinizations were classified in \cite{Chari2, Cha7, Cha8, Cha9} for all types, except for type $D,E$ for a weight orthogonal to the special node. For the regular cases (ie. with a linear Dynkin diagram, that is to say types $A$, $B$, $C$, $F_4$, $G_2$), the classification is complete :

\begin{thm}\cite{Chari2, Cha7, Cha8}\label{class} Suppose that $\Glie$ is regular and let $\lambda\in P^+$. For $i\in I$ let $\lambda_i=\lambda(\alpha_i^{\vee})$ and for $i < n$ let
$$c_i(\lambda)=q^{r_i\lambda_i+r_{i+1}\lambda_{i+1}+r_{i+1}-C_{i+1,i}-1}\text{ and }c_i'(\lambda)=q^{r_i\lambda_i+r_{i+1}\lambda_{i+1}+r_{i}-C_{i,i+1}-1}.$$ 
Then a simple $\U_q(\Lo\Glie)$-module $L(m)$ is a minimal affinization of $V(\lambda)$ if and only if $m$ is of the form $m = \prod_{i\in I} X_{\lambda_i,a_i}^{(i)}$ with $(a_i)_{i\in I}\in(\CC^*)^I$ satisfying one of two conditions :

(I) For all $i<j\in I$, $a_i/a_j=\prod_{i\leq s <j}c_s(\lambda)$.

(II) For all $i<j\in I$, $a_j/a_i=\prod_{i\leq s <j}c_s'(\lambda)$.\end{thm}

\noindent Observe that we have rewritten the defining formulas of $c_s$, $c_s'$ \cite{Chari2, Cha7, Cha8} in a slightly different (more homogeneous) way.

\noindent Observe that for classical types, minimal affinizations (called generalized Kirillov-Reshetikhin modules) were also studied in \cite{GK}.

\begin{rem}\label{kraff} As a consequence of Theorem \ref{class}, for $k\geq 0$ and $i\in I$, the minimal affinizations of $V(k\Lambda_i)$ are the Kirillov-Reshetikhin modules.\end{rem} 

For $\Glie$ of type $D$, and $\lambda\in P^+$, we define with the same formulas $c_i(\lambda)$ for $i < n-1$, and we set $c_{n-1}(\lambda)=q^{\lambda_{n-2}+\lambda_n + 1}$. For a monomial $m = \prod_{i\in I} X_{\lambda_i,a_i}^{(i)}$ we have analog conditions (I) and (II) :

(I) For all $i<j\in I$, $a_i/a_j=(c_{j-1}(\lambda))^{\epsilon_j}\prod_{i\leq s \leq \text{min}(j-1,n-3)}c_s(\lambda)$,

(II) For all $i<j\in I$, $a_j/a_i=(c_{j-1}(\lambda))^{\epsilon_j}\prod_{i\leq s \leq \text{min}(j-1,n-3)}c_s(\lambda)$,
\\where $\epsilon_j = 0$ if $j\leq n-2$ and $\epsilon_{n-1} = \epsilon_n = 1$.

\noindent It follows from \cite[Theorem 6.1]{Cha7} that if $\lambda_{n-2}\neq 0$ and $m = \prod_{i\in I} X_{\lambda_i,a_i}^{(i)}$ satisfies $(I)$ or $(II)$, then $L(m)$ is a minimal affinization of $V(\lambda)$.

\subsection{Main results}

It follows directly from Theorem \ref{formerkr} and remark \ref{kraff} that (see also Proposition \ref{qseg} for an alternative general proof) :

\begin{cor} For $i\in I$ and $k\geq 0$, the minimal affinizations of $V(k\Lambda_i)$ are special.\end{cor}

In general a minimal affinization is not special. Let us look at some examples.

First we consider type $C$. 

If $m$ satisfies condition $(II)$ of Theorem \ref{class}, $L(m)$ is not necessarily special. For example consider the case $\Glie$ of type $C_3$ and $m=Y_{2,1}Y_{2,q^2}Y_{3,q^7}$. $L(m)$ is a minimal affinization of $V(2\Lambda_2+\Lambda_3)$. By the process described in remark \ref{process}, the monomials $1_1 1_3 2_2^{-1} 2_4^{-1} 3_1 3_3 3_7$, $1_3^{-1} 1_5^{-1} 3_1 3_3 3_7$, $1_3^{-1} 1_5^{-1} 2_2 2_4 3_5^{-1} 3_3 3_7$, $2_6^{-1} 2_4^{-1} 3_3^2 3_7$ and $3_3$ occur in $\chi_q(L(m))$ and so $L(m)$ is not special.

If $m$ satisfies condition $(I)$ of Theorem \ref{class}, $L(m)$ is not necessarily special. For example consider the case $\Glie$ of type $C_3$ and $m=Y_{1,q^3}Y_{1,q^5}Y_{1,q^7}Y_{2,1}$. $L(m)$ is a minimal affinization of $V(3\Lambda_1+\Lambda_2)$. By the process described in remark \ref{process}, the monomials $1_1 1_3 1_5 1_7 2_2^{-1} 3_1$, $1_1 1_3 1_5 1_7 2_4 3_5^{-1}$, $1_1 1_3 1_5^2 1_7 2_6^{-1}$, $1_1 1_3 1_5$ occur in $\chi_q(L(m))$ and so $L(m)$ is not special. 

Eventually, let $m = Y_{1,1}Y_{1,q^2}Y_{1,q^4}Y_{2,q^7}Y_{2,q^9}Y_{3,q^{14}}$. We can see as for $Y_{2,1}Y_{2,q^2}Y_{3,q^7}$ that $L(m)$ is not special. Moreover $L(m)$ is antispecial if and only if the module $L(Y_{1,q^{14}}Y_{1,q^{12}}Y_{1,q^{10}}Y_{2,q^7}Y_{2,q^5}Y_{3,1})$ is special (see Lemma \ref{involchiq} and Corollary \ref{invmon} bellow). But we can check as for $Y_{1,q^3}Y_{1,q^5}Y_{1,q^7}Y_{2,1}$ that this module is not special. So $L(m)$ is not special and not antispecial.

For types $D$, there are minimal affinizations which are not special. For example let $\Glie$ of type $D_4$ and $m = Y_{1,q^3}Y_{1,q^5}Y_{2,1}$. Then $L(m)$ is not special (see \cite[Remark 6.8]{small}).

\noindent However we prove in this paper :

\begin{thm}\label{fora} For $\Glie$ of type $A$, $B$, $G$, all minimal affinizations are special and antispecial.\end{thm}

\begin{thm}\label{forcd} For $\Glie$ of type $C$, $F_4$ and $\lambda\in P$ satisfying $\lambda_n = 0$, all minimal affinizations of $V(\lambda)$ satisfying (I) (resp. (II)) are antispecial (resp. special).

For $\Glie$ of type $D$ and $\lambda\in P$ satisfying $\lambda_{n-1} = \lambda_n$, all $L(m)$ satisfying (I) (resp. (II)) and $a_{n-1} = a_n$ are antispecial (resp. special).\end{thm}

\noindent Note that for type $D$, the condition $a_{n-1} = a_n$ is automatically satisfied if $\lambda_j\neq 0$ for one $j\leq n-2$.

\begin{thm}\label{egalun} For $\Glie$ of type $A$, $B$, all minimal affinizations are thin.\end{thm}

\noindent Theorems \ref{fora}, \ref{forcd} and \ref{egalun} are proved in section \ref{proofmainres}.

Note for type $C$, there are minimal affinizations which are not thin : for example consider $\Glie$ of
type $C_4$ and $m=Y_{3,1}Y_{3,q^2}$. $L(m)$ is a Kirillov-Reshetikhin and a minimal affinization of $V(2\Lambda_3)$. By the process described in remark \ref{process}, the following monomials occur in $\chi_q(L(m))$ : $3_1 3_2$, $2_1 2_3 3_2^{-1} 3_4^{-1} 4_1 4_3$, $2_1 2_3 4_5^{-1} 4_3$, $1_4 2_1 2_5^{-1} 3_4 4_5^{-1} 4_3$, 
$1_4 2_1 3_6^{-1} 4_3$, 
$1_4 2_1 3_4 4_7^{-1}$, 
$1_6^{-1} 2_1 2_5 3_4 4_7^{-1}$ and $1_6^{-1} 2_1 2_5^2 3_6 4_5 4_7^{-1}$. And so by Proposition \ref{jdecomp} and Proposition \ref{aidesldeux} the monomial $2_1 2_5 2_7^{-1} 3_6^2 4_5 4_7^{-1}$ occurs in $\chi_q(L(m))$ with multiplicity larger than $2$.

For type $G_2$, there are minimal affinizations which are not thin : for example let $m=Y_{2,0}Y_{2,2}$. $L(m)$ is a Kirillov-Reshetikhin and a minimal affinization of $V(2\Lambda_2)$. We have 
$2_0 2_2$, $1_1 1_3 2_4^{-1} 2_2^{-1}$, $1_7^{-1}1_9^{-1}2_42_6^22_8\in\mathcal{M}(L(m))$,
and so $Y_{1,9}^{-1}Y_{2,4}Y_{2,6}$ occurs in $\chi_q(L(m))$ with multiplicity larger than $2$.

\section{Preliminary results}\label{trois}

In this section $\Glie$ is an arbitrary semi-simple Lie algebra. We discuss preliminary results which will be used in the proof of Theorem \ref{fora}, \ref{forcd} and \ref{egalun} in the next section.

First it is well known that :

\begin{lem}\label{produit} Let $L(m_1), L(m_2)$ be two simple modules. Then $L(m_1m_2)$ is a subquotient of $L(m_1)\otimes L(m_2)$. In particular $\mathcal{M}(L(m_1m_2))\subset\mathcal{M}(L(m_1))\mathcal{M}(L(m_2))$.\end{lem}

\subsection{Results of \cite{small}}

All results of this subsection are preliminary results of \cite{small}.

\begin{lem}\label{plusgrand} Let $a\in\CC^*$ and $m$ be a monomial of $\ZZ[Y_{i,aq^r}]_{i\in I,r\geq 0}$. Then for $m'\in\mathcal{M}(L(m))$ and $b\in\CC^*$, ($v_{i,b}(m'm^{-1})\neq 0\Rightarrow b\in aq^{r_i + \NN}$).\end{lem}

\begin{lem}\label{thinmon} Let $V\in\text{Rep}(\U_q(\Lo\Glie))$ be a $\U_q(\Lo\Glie)$-module and $m'\in\mathcal{M}(V)$ such that there is $i\in I$ satisfying $\text{Min}\{u_{i,a}(m')|a\in\CC^*\}\leq -2$. Then there is $m''>m'$ in $\mathcal{M}(V)$ $i$-dominant such that $\text{Max}\{u_{i,b}(m'')|b\in\CC^*\} \geq 2$.\end{lem}

We recall \cite{small} that a monomial $m$ is said to be thin if $\text{Max}_{i\in I, a\in\CC^*} |u_{i,a}(m)|\leq 1$.

\begin{prop}\label{proofthin} If $V$ is thin then all $m\in\mathcal{M}(V)$ are thin. If $V$ is special and all $m\in\mathcal{M}(V)$ are thin, then $V$ is thin.\end{prop}

\begin{lem}\label{depart} Let $m$ dominant and $J\subset I$. Let $v$ be an highest weight vector of $L(m)$ and $L'\subset L(m)$ the $\U_q(\Lo \Glie_J)$-submodule of $L(m)$ generated by $v$. Then $L'$ is an $\U_q(\Lo\Hlie)$-submodule of $L(m)$ and $\chi_q(L')=L_J(m)$.\end{lem}

\noindent In particular for $\mu\in \omega(m) - {\sum}_{j\in J}\NN\alpha_j$, we have $$\text{dim}((L(m))_\mu)=\text{dim}((L^J(m^{\rightarrow (J)}))_{\mu^{\rightarrow (J)}}),$$ 
where $\mu^{\rightarrow (J)} = {\sum}_{j\in J}\mu(\alpha_j^{\vee})\omega_j$.

\begin{lem}\label{oudeux} Let $V = L(m)$ be a $\U_q(\Lo \Glie)$-module simple module and consider a monomial $m'\in(\mathcal{M}(L(m))-\{m\})$. Then there is $j\in I$ and $M'\in\mathcal{M}(V)$ $j$-dominant such that $M'>m'$, $M'\in m'\ZZ[A_{j,b}]_{b\in\CC^*}$ and $((\U_q(\Lo\Glie_j).V_{M'})\cap (M)_{m'})\neq \{0\}$.\end{lem}

We have the following elimination theorem :

\begin{thm}\label{racourc} Let $V = L(m)$ be a $\U_q(\Lo \Glie)$-module simple module. Let $m' < m$ and $i\in I$ satisfying the following conditions 

(i) there is a unique $i$-dominant $M\in(\mathcal{M}(V)\cap m'\ZZ[A_{i,a}]_{a\in\CC^*})-\{m'\}$,

(ii) $\sum_{r\in\ZZ} x_{i,r}^+ (V_M)=\{0\}$,

(iii) $m'$ is not a monomial of $L_i(M)$,

(iv) if $m''\in\mathcal{M}(\U_q(\Lo\Glie_i).V_M)$ is $i$-dominant, then $v(m''m^{-1}) \geq v(m'm^{-1})$,

(v) for all $j\neq i$, $\{m''\in\mathcal{M}(V)|v(m''m^{-1}) < v(m'm^{-1})\}\cap m'\ZZ[A_{j,a}^{\pm 1}]_{a\in\CC^*}=\emptyset$.

\noindent Then $m'\notin \mathcal{M}(V)$.\end{thm}

\begin{lem}\label{etoile}\label{stara}\label{remont} Let $L(m)$ be a simple $\U_q(\Lo\Glie)$-module, and $m'\in\mathcal{M}(L(m))$ such that all $m''\in\mathcal{M}(L(m))$ satisfying $v(m''m^{-1}) < v(m'm^{-1})$ is thin. 

1) For $i\in I$ such that $m'$ is not $i$-dominant, there is $a\in\CC^*$ such that $u_{i,a}(m') < 0$ and $m'A_{i,aq_i^{-1}}\in\mathcal{M}(L(m))$.

2) We suppose that $\Glie=sl_{n+1}$, that there are $i\in I$, $a\in\CC^*$ satisfying $u_{i,a}(m') = -1$ and $m'Y_{i,a}$ is dominant. Then there is $M\in\mathcal{M}(L(m))$ dominant such that $M > m'$ and $v_n(m'M^{-1})\leq 1$, $v_1(m'M^{-1})\leq 1$.

3) We suppose that $\Glie=sl_{n+1}$, that there is $j\in I$, such that $m'$ is $(I-\{j\})$-dominant and if $j\leq n-1$, then for all $a\in\CC^*$, $(u_{j,a}(m')<0\Rightarrow u_{j+1,aq^{-1}}(m') > 0)$. Then there is $M\in\mathcal{M}(L(m))$ dominant of the form 
$$M = m'\prod_{\{a\in\CC^*|u_{j,a}(m')<0\}}(A_{j,aq^{-1}}A_{j-1,aq^{-3}}\cdots A_{i_a,aq^{i_a-j-1}}),$$
where for $a\in\CC^*$, $1\leq i_a\leq j$.
\end{lem}

\noindent Observe that we can prove in the same way an analog result where we replace all $i\in I$ by $\overline{i}=n-i+1$.

\subsection{Involution of $\U_q(\Lo\Glie)$ and simple modules}\label{secinvol} For $\mu$ an automorphism of $\U_q(\Lo\Glie)$ and $V$ a $\U_q(\Lo\Glie)$-module we denote the corresponding twisted module by $\mu^*V$. The involution of the algebra $\Yim$ defined by $Y_{i,a}^{\pm}\mapsto Y_{i,a^{-1}}^{\mp}$ is denoted by $\sigma$. 

For all $b\in\CC^*$, let $\tau_b$ be the automorphism of $\U_q(\Lo\Glie)$ defined by $x_{i,m}^{\pm}\mapsto b^{-m}x_{i,m}^{\pm}$, $h_{i,r}\mapsto b^{-r} h_{i,r}$, $k_i^{\pm}\mapsto k_i^{\pm}$. For $V$ a $\U_q(\Lo\Glie)$-module we have $\chi_q(\tau_b^*V)=\beta_b(\chi_q(V))$ where $\beta_b:\Yim\rightarrow \Yim$ is the ring morphism such that $\beta_b(Y_{i,a}^{\pm})=Y_{i,ab}^{\pm}$. So $\tau_b^*L(m)\simeq L(\beta_b(m))$ and $\chi_q(\tau_b^*L(m))=\beta_b(\chi_q(L(m)))$. 

\begin{lem}\cite[Proposition 1.6]{Chari2} There is a unique involution $\sigma$ of the algebra $\U_q(\Lo\Glie)$ such that for all $i\in I, r\in\ZZ, m\in\ZZ-\{0\}$ : 
$$\sigma(x_{i,r}^{\pm})=x_{i,-r}^{\mp}\text{ , }\sigma(h_{i,m})=-h_{i,-m}\text{ , }\sigma(k_i)=k_i^{-1}.$$ 
Moreover for $m\geq 0$, $\sigma(\phi_{i,\pm m}^{\pm})=\phi_{i,\mp m}^{\mp}$.
\end{lem}

(Observe that we could also use $\sigma(x_{i,r}^{\pm})=-x_{i,-r}^{\mp}$ to define an involution of $\U_q(\Lo\Glie)$.)

\begin{lem}\label{involchiq} We have $\chi_q(\sigma^*V) = \sigma(\chi_q(V))$.\end{lem}

\demo For $\gamma = (\gamma_{i,\pm m}^{\pm})_{i\in I,m\geq 0}$, it follows from the relation $\sigma(\phi_{i,\pm m}^{\pm})=\phi_{i,\mp m}^{\mp}$ that $V_{\gamma} = (\sigma^* V)_{\gamma'}$ where $\gamma'=(\gamma_{i,\mp m}^{\pm})_{i\in I,m\geq 0}$. Let $Q_i(u)=\prod_{a\in\CC^*}(1-ua)^{q_{i,a}}$ and $R_i(u)=\prod_{a\in\CC^*}(1-ua)^{r_{i,a}}$ such that in $\CC[[u^{\pm}]]$ we have :
\begin{equation*}
{\sum}_{m\geq 0}\gamma_{i,\pm m}^{\pm}u^{\pm 
m}=q_i^{\text{deg}(Q_i)-\text{deg}(R_i)}\frac{Q_i(uq_i^{-1})R_i(uq_i)}{Q_i(uq_i)R_i(uq_i^{-1})}.
\end{equation*}
Then in $\CC[[u^{\pm}]]$ we have :
\begin{equation*}
\begin{split}
\sum_{m\geq 0}\gamma_{i,\mp m}^{\mp}u^{\pm 
m}&=q_i^{\text{deg}(Q_i)-\text{deg}(R_i)}\frac{Q_i(u^{-1}q_i^{-1})R_i(u^{-1}q_i)}{Q_i(u^{-1}q_i)R_i(u^{-1}q_i^{-1})}\\&=q_i^{\text{deg}(Q_i')-\text{deg}(R_i')}\frac{Q_i'(uq_i^{-1})R_i'(uq_i)}{Q_i'(uq_i)R_i'(uq_i^{-1})},
\end{split}
\end{equation*}
where $Q_i'(u)=\prod_{a\in\CC^*}(1-ua)^{r_{i,a^{-1}}}$ and $R_i'(u)=\prod_{a\in\CC^*}(1-ua)^{q_{i,a^{-1}}}$ by using the identities
\begin{equation*}
q_i\frac{1-au^{-1}q_i^{-1}}{1-au^{-1}q_i}=q_i^{-1}\frac{1-a^{-1}uq_i}{1-a^{-1}uq_i^{-1}}\text{ and }
q_i^{-1}\frac{1-au^{-1}q_i}{1-au^{-1}q_i^{-1}}=q_i\frac{1-a^{-1}uq_i^{-1}}{1-a^{-1}uq_i}.
\end{equation*}
\qed

In particular $\chi(\sigma^*V)=\sigma(\chi(V))$ where $\sigma:\mathcal{E}\rightarrow\mathcal{E}$ is defined by $\sigma(e(\lambda))=e(-\lambda)$.

Let $w_0$ be the longest element in the Weyl group of $\Glie$ and $i\mapsto \overline{i}$ be the unique bijection of $I$ such that $w_0(\alpha_i) = -\alpha_{\overline{i}}$. Let $h^{\vee}$ be the dual Coxeter number of $\Glie$ and $r^{\vee}$ the maximal number of edges connecting two vertices of the Dynkin diagram of $\Glie$.

\begin{cor}\label{invmon} For $m$ dominant, we have $\sigma^*L(m)\simeq L(m')$ where
$$m'=\prod_{a\in\CC^*}{\prod}Y_{\overline{i},a^{-1}q^{-r^{\vee}h^{\vee}}}^{u_{i,a}(m)}.$$
\end{cor}

\demo A submodule of $V$ is a submodule of $\sigma^*V$, so $V$ simple implies $\sigma^*V$ simple. As it is proved in \cite[Corollary 6.9]{Fre2} that the lowest monomial of $L(m)$ is $\prod_{i\in I, a\in\CC^*}Y_{\overline{i},aq^{r^{\vee}h^{\vee}}}^{-u_{i,a}(m)}$, we get the result from Lemma \ref{involchiq}.\qed

Observe that as a by product we get the following symmetry property : 

\begin{cor} If $k=\overline{k}$, then $\chi_q(W_{k,a}^{(i)})$ is invariant by $(\tau_{a^2q^{r^{\vee}h^{\vee}}}\circ\sigma)$.\end{cor}

\noindent For example, this symmetry can be observe on the diagrams of $q$-characters in \cite[Figure 1]{Naex} and \cite[Section 8]{her02}.

Let us go back to the main purposes of this paper. First we get a simplification in the proof of Theorem \ref{fora} :

\begin{cor}\label{simplification} In Theorem \ref{fora}, it suffices to prove that all minimal affinizations are special.\end{cor}

\demo First suppose that $\Glie$ is of type $B$ or $G$. Then $\overline{i}=i$. If $m$ satisfies condition $(II)$ of Theorem \ref{class}, then $m'$ of corollary \ref{invmon} satisfies condition $(I)$. Moreover if $M$ is dominant, then $\sigma(M)$ is antidominant. So we can conclude with Lemma \ref{involchiq}. If $\Glie$ is of type $A$, conditions (I) and (II) are the same up to a different numbering.\qed

Exactly in the same way we get :

\begin{cor}\label{simplificationcd} In Theorem \ref{forcd}, it suffices to prove that the considered simple representations satisfying the condition (II) are special.\end{cor}

\noindent For $V$ a $\U_q(\Lo\Glie)$-module, denote by $V^*$ the dual module of $V$. As $S(k_i)=k_i^{-1}$, we have $\chi(V^*)=\sigma(\chi(V))$. As a direct consequence of \cite[Corollary 6.9]{Fre2}, we have :

\begin{lem} For $m$ dominant, we have $(L(m))^*\simeq L(m'')$ where 
$$m'' = {\prod}_{i\in I,a\in\CC^*}Y_{\overline{i},aq^{-r^{\vee}h^{\vee}}}^{u_{i,a}(m)}.$$
\end{lem}

Note that it was proved in \cite{Fre2} that we have the following relation between the $q$-character of $(V_i(a))^*\simeq V_{\overline{i}}(aq^{-r^{\vee}h^{\vee}})$ and $V_i(a)$ : $$\chi_q((V_i(a))^*)=(\tau_a\circ\sigma\circ\tau_{a^{-1}})(\chi_q(V_i(a))).$$

\begin{prop}\label{vraiinv} For $m$ a dominant monomial, we have 
$$\chi(L(m))=\chi(L((\sigma(m))^{-1})).$$\end{prop}

\demo From previous results, we have 
$$\chi(\sigma^*((L(m))^*))=\sigma(\chi((L(m))^*))=\chi(L(m)),$$ 
and $\sigma^*((L(m))^*)\simeq L({\prod}_{i\in I, a\in\CC^*}Y_{i,a^{-1}}^{u_{i,a}(m)}) = L((\sigma(m))^{-1})$.\qed

The above proposition can be extended to $\chi(L(m)) = \chi(L({\prod}_{i\in I, a\in\CC^*}Y_{i,ba^{-1}}^{u_{i,a}(m)}))$ for all $b\in\CC^*$. 

Observe that we do not have a direct relation between the monomials of the same weight space : for example for $\Glie = sl_2$ and $m=Y_qY_{q^3}^2$, the term of weight $\Lambda$ in $\chi_q(L(m))$ is $2Y_qY_{q^3}Y_{q^5}^{-1}$ and the term of weight $\Lambda$ in $\chi_q(L(\sigma(m)))$ is $Y_{q^{-3}}+Y_{q^{-3}}^2Y_q^{-1}$.

\subsection{Additional preliminary results}

\begin{lem}\label{geneun} Let $m=X_{k,a}^{(i)}$. Let $m'\in \mathcal{M}(W_{k,a}^{(i)})$ and $\mu\in\{k,k - 2,\cdots,-k+2\}$. Then $v_{i,aq_i^{\mu}}(m'm^{-1})\geq 1$ implies 
\begin{equation*}
v_{i,aq_i^k}(m'm^{-1})\geq 1\text{ , } v_{i,aq_i^{k-2}}(m'm^{-1})\geq 1\text{ , } \cdots \text{ , } v_{i,aq_i^{\mu}}(m'm^{-1})\geq 1.
\end{equation*}
\end{lem}

\demo For $\mu=k$ the result is clear. We suppose that $\mu < k$ and we prove the result by induction on $k$. For $k=1$ the result is clear. For general $k\geq 1$ and $\mu < k$, suppose that $v_{i,aq_i^{\mu}}(m'm^{-1})\geq 1$. So $m'\neq m$ and it follows from Theorem \ref{formerkr} that $m'\leq mA_{i,aq_i^k}^{-1}$. By Lemma \ref{produit}, $m'=m_1m_2$ where $m_1\in\mathcal{M}(V_i(aq_i^{k-1}))$ and $m_2\in \mathcal{M}(W_{k-1,aq_i^{-1}}^{(i)})$. From Lemma \ref{plusgrand}, $v_{j,b}(m_1Y_{i,aq_i^{k-1}}^{-1})\neq 0$ implies $b = aq^{r_i(k-1)+R}$ with $R\geq 1$ and so $b = aq_i^{\mu}$. So we have $v_{i,aq_i^{\mu}}(m_2(X_{k-1,aq_i^{-1}}^{(i)})^{-1})\geq 1$. So by the induction hypothesis 
\begin{equation*}
\begin{split}
v_{i,aq_i^{k-2}}(m_2(X_{k-1,aq_i^{-1}}^{(i)})^{-1})\geq 1\text{ , } &v_{i,aq_i^{k-4}}(m_2(X_{k-1,aq_i^{-1}}^{(i)})^{-1})\geq 1,
\\&\cdots, v_{i,aq_i^{\mu}}(m_2(X_{k-1,aq_i^{-1}}^{(i)})^{-1})\geq 1.
\end{split}
\end{equation*} 
We can conclude because it follows from Theorem \ref{formerkr} that $v_{i,aq_i^k}(m'm^{-1})\geq 1$.\qed

\begin{lem}\label{plushaut} Let $a\in\CC^*$ and a monomial $m\in\ZZ[Y_{i,aq^r}]_{i\in I,r\in\ZZ}$. Let $m'\in\mathcal{M}(L(m))$ and $R\in\ZZ$ such that for all $i\in I$, $(u_{i,aq^r}(m') < 0\Rightarrow r\leq R)$. Then there is a dominant monomial $M\in\mathcal{M}(L(m))\cap (m\ZZ[A_{i,aq^r}]_{\{(i,r)|i\in I, r\leq R-r_i\}})$. \end{lem}

\demo From Lemma \ref{jdecomp} it suffices to prove the result for $\U_q(\hat{sl_2})$. In this case the result follows from (3) of Proposition \ref{aidesldeux}.\qed

To conclude this section, let us prove a refined version of Proposition \ref{jdecomp}. For $i\in I$, $a\in\CC^*$ and $m$ a monomial denote $m^{\rightarrow (i,a)}={\prod}_{r\in\ZZ}Y_{i,aq_i^{2r}}^{u_{i,aq_i^{2r}}(m)}$.  Define :
$$L_{i,a}(m)=m \mu_i^I((m^{\rightarrow (i,a)})^{-1}\chi_q^i(L^i(m^{\rightarrow (i,a)}))).$$
Observe that for $a'\in aq_i^{2\ZZ}$, $m^{\rightarrow (i,a)}=m^{\rightarrow (i,a')}$ and $L_{i,a}(m)=L_{i,a'}(m)$. So the definition can be given for $a\in(\CC^*/q_i^{2\ZZ})$. We have :

\begin{cor}\label{jdecompprime} For a representation $V\in\text{Rep}(\U_q(\Lo \Glie))$, $i\in I$ and $a\in\CC^*$, there is a unique decomposition in a finite sum : 
$$\chi_q(V)={\sum}_{\{m'|(m')^{\rightarrow (i,a)}\text{ is dominant}\}}\lambda_{i,a}(m')L_{i,a}(m').$$
Moreover for all $m'$ such that $(m')^{\rightarrow (i,a)}$ is dominant, we have $\lambda_{i,a}(m')\geq 0$.
\end{cor}

\demo First we write the decomposition of Lemma \ref{jdecomp} with $J=\{i\}$. Then it follows from Proposition \ref{aidesldeux} that for $m'$ an $i$-dominant monomial we have 
$$L_i(m')=(m')^{\rightarrow (I-\{i\})}\prod_{b\in (\CC^*/q_i^{2\ZZ})}L_{i,b}((m')^{\rightarrow (i,b)}).$$\qed

\section{Proof of the main results}\label{proofmainres}

In this section we prove Theorems \ref{fora}, \ref{forcd} and \ref{egalun}. We study successively the different types.

\subsection{Type $A$}\label{typa} 

In this section \ref{typa}, $\Glie=sl_{n+1}$. 

\begin{lem}\label{proofspe} Let $\lambda\in P^+$ and $L(m)$ be a minimal affinization of $V(\lambda)$. Suppose that $m$ satisfies the condition $(II)$ (resp. condition $(I)$) of Theorem \ref{class}. Let $K=\text{max}\{i\in I|\lambda_i\neq 0\}$ (resp. $K=\text{min}\{i\in I|\lambda_i\neq 0\}$). The following properties are satisfied.
\begin{enumerate}
\item For all $m'\in\mathcal{M}(L(m))$, if $v_K(m'm^{-1})\geq 1$, then $v_{K,a_Kq^{\lambda_K}}(m'm^{-1})\geq 1$.

\item $L(m)$ is special.

\item $L(m)$ is thin.

\item For all $m'\in \mathcal{M}(L(m))$ we have 
\begin{equation*}
\begin{split}
v_{j,a_kq^{\lambda_k+|j-k|}}(m'm^{-1})&=v_{j,a_kq^{\lambda_k+|j-k|-2}}(m'm^{-1})
\\&=\cdots=v_{j,a_kq^{\lambda_k+|j-k|-2R}}(m'm^{-1})=1,
\end{split}
\end{equation*}
where 
$$j=\text{max}\{i|v_i(m'm^{-1})\neq 0\}\text{ (resp. $j=\text{min}\{i|v_i(m'm^{-1})\neq 0\}$)},$$ 
$$k=\text{max}\{i\leq j|\lambda_i\neq 0\}\text{ (resp. $k=\text{min}\{i\geq j|\lambda_i\neq 0\}$)},$$ 
and $R=v_j(m'm^{-1})-1$.
\end{enumerate}
\end{lem}

Observe that as a consequence of property (4), for $b\in\CC^*$, $v_{j,b}(m'm^{-1})\neq 0$ implies 
\begin{equation*}
b\in\{a_kq^{\lambda_k+|j-k|},a_kq^{\lambda_k+|j-k|-2},\cdots,a_kq^{\lambda_k+|j-k|-2R}\}.
\end{equation*}

Lemma \ref{proofspe} combined with corollary \ref{simplification} implies Theorem \ref{fora} and Theorem \ref{egalun} for type $A$.

\demo We suppose that $L(m)$ satisfies $(II)$ (the case $(I)$ is treated in the same way). We prove by induction on $u(m)\geq 0$ simultaneously that (1), (2), (3) and (4) are satisfied.

For $u(m)=0$ the result is clear. Suppose that $u(m)\geq 1$. 

First we prove (1) by induction on $v(m'm^{-1})\geq 0$. For $v(m'm^{-1})=0$ we have $m'=m$ and the result is clear. In general suppose that for $m''$ such that $v(m''m^{-1}) < v(m'm^{-1})$, the property is satisfied. Suppose that $v_K(m'm^{-1})\geq 1$ and $v_{K,a_Kq^{\lambda_K}}(m'm^{-1}) = 0$. It suffices to prove that the conditions of Proposition \ref{racourc} with $i = K$ are satisfied. 

Condition (i) of Proposition \ref{racourc} : if $M>m'$ and $M\in\mathcal{M}(L(m))$, we have $v_{K,a_Kq^{\lambda_K}}(Mm^{-1}) = 0$ and so by the induction hypothesis $v_{K}(Mm^{-1}) = 0$. So if we suppose moreover that $M\in m'\ZZ[A_{K,a}]_{a\in\CC^*}$, we have $M=m'\prod_{a\in\CC^*}A_{K,a}^{v_{K,a}(m'm^{-1})}$, and so we get the uniqueness. For the existence, it suffices to prove that this $M=m'\prod_{a\in\CC^*}A_{K,a}^{v_{K,a}(m'm^{-1})}$ is in $\mathcal{M}(L(m))$. By Lemma \ref{oudeux}, there is $j\in I$, $M'\in\mathcal{M}(L(m))$ $j$-dominant such that $M'>m'$ and $M'\in m'\ZZ[A_{j,a}]_{a\in\CC^*}$. By the induction hypothesis on $v$ we have $j=K$, and so by uniqueness $M'=M$. 

Condition (ii) of Proposition \ref{racourc} : by construction of $M$ we have $v_K(Mm^{-1})=0$. 

Condition (iii) of Proposition \ref{racourc} : first observe that 
$$M\in m^{\rightarrow (K)}\mathcal{M}(L(m(m^{\rightarrow (K)})^{-1})).$$ 
As $u(m(m^{\rightarrow (K)})^{-1}) < u(m)$, we have property (4) for $L(m(m^{\rightarrow (K)})^{-1})$ and we get 
$$(M)^{\rightarrow (K)} = Y_{K,a_K q^{\lambda_K-1}}Y_{K,a_K q^{\lambda_K-3}}\cdots Y_{K,a_K q^{-\lambda_K+1-2R'}},$$ 
with $R'\geq 0$. By Lemma \ref{aidesldeux}, $m'$ is not a monomial of $\mathcal{M}(L_K(M))$. 

Condition (iv) of Proposition \ref{racourc} : let $m''\in\mathcal{M}(\U_q(\Lo\Glie_K).(L(m))_M)$ such that $v(m''m^{-1}) < v(m'm^{-1})$. Then we have $m''\in MA_{K,a_Kq^{\lambda_k}}^{-1}\ZZ[A_{K,b}^{-1}]_{b\in\CC^*}$ and so $(m'')^{\rightarrow (K)}$ is right negative, so $m''$ is not $K$-dominant. 

Condition (v) of Proposition \ref{racourc} : clear by the induction property on $v$.

Now we prove (2). Let $J=\{i\in I|i<K\}$. By Lemma \ref{produit}, $\mathcal{M}(L(m))\subset (m^{\rightarrow (J)}\mathcal{M}(L(m^{\rightarrow (K)})))\cup (\mathcal{M}(L(m^{\rightarrow (J)}))m^{\rightarrow (K)})$. From Theorem \ref{formerkr}, all monomials of $m^{\rightarrow (J)}(\chi_q(L(m^{\rightarrow (K)}))-m^{\rightarrow (K)})$ are lower than $mA_{K,a_Kq^{\lambda_K}}^{-1}$ which is right-negative, and so are not dominant. Let $m'\in(\mathcal{M}(L(m^{\rightarrow (J)}))m^{\rightarrow (K)}-\{m\})$. If $v_K(m'm^{-1})\geq 1$, it follows from property (1) that $m'$ is lower than $mA_{K,a_Kq^{\lambda_K}}^{-1}$ which is right-negative, so $m'$ is not dominant. If $v_K(m'm^{-1})=0$, we have $u_{K,b}(m'(m^{\rightarrow (K)})^{-1})\geq 0$ for all $b\in\CC^*$. We have $m'(m^{\rightarrow (K)})^{-1}\in\mathcal{M}(L(m^{\rightarrow (J)}))$ with $u(m^{\rightarrow (J)}) < u(m)$, so by the induction hypothesis on $u$, $m'(m^{\rightarrow (K)})^{-1}$ is not dominant. So there is $i\neq K$, $b\in\CC^*$, such that $u_{i,b}(m'(m^{\rightarrow (K)})^{-1}) < 0$. As $u_{i,b}(m'(m^{\rightarrow (K)})^{-1})=u_{i,b}(m')$, $m'$ is not dominant. So $L(m)$ is special.

Now we prove (3). From property (2) and Proposition \ref{proofthin}, it suffices to prove that all monomials of $\mathcal{M}(L(m))$ are thin. From Lemma \ref{thinmon}, we can suppose that there is $m'\in\mathcal{M}(L(m))$ such that there are $i\in I, a\in\CC^*$ satisfying $u_{i,a}(m')=2$ and such that all $m''$ satisfying $v(m''m^{-1}) < v(m'm^{-1})$ is thin. Then $m'$ is $(\{1,\cdots,i-2\}\cup \{i\}\cup \{i+2,\cdots,n\})$-dominant and $(u_{i-1,b}(m')<0\Rightarrow b=aq)$ and $(u_{i+1,b}(m')<0\Rightarrow b=aq)$. We can apply (3) of Lemma \ref{remont} for $\Glie_{\{1,\cdots,i-1\}}$ and for $\Glie_{\{i+1,\cdots,n\}}$. We get $M\in\mathcal{M}(L(m))$ dominant satisfying $u_{j_1,aq^{j_1-i}}(M)\geq 1$, $u_{j_2,aq^{i-j_2}}(M)\geq 1$ with $j_1 < j_2$, $j_1\leq i\leq j_2$. From property (2) we have $m=M$, contradiction with condition (II) of Theorem \ref{class}. 

Now we prove (4) by induction on $v(m'm^{-1})\geq 0$. We can suppose that $j=n$ (Lemma \ref{depart}). So $k=K$. For $v(m'm^{-1})=0$ we have $m'=m$ and the result is clear. Let $m'$ such that the property is satisfied for $m''$ with $v(m''m^{-1})<v(m'm^{-1})$. Let $R\geq 0$ maximal such that 
$$m'\leq m A_{n,a_kq^{\lambda_k+n-k}}^{-1}A_{n,a_kq^{\lambda_k+n-k-2}}^{-1}\cdots A_{n,a_kq^{\lambda_k+n-k-2R+2}}^{-1}.$$ 
We suppose moreover that 
$$m'\leq m A_{n,a_kq^{\lambda_k+n-k}}^{-1}A_{n,a_kq^{\lambda_k+n-k-2}}^{-1}\cdots A_{n,a_kq^{\lambda_k+n-k-2R+2}}^{-1}A_{n,b}^{-1}$$ 
with $b\neq a_kq^{\lambda_k+n-k-2R}$. 
By the induction property on $v$, $m'$ is $(I-\{n\})$-dominant, $u_{n,bq}(m') < 0$ and $(u_{n,c}(m')<0\Rightarrow c=bq)$. By property (3), $u_{n,bq}(m') = -1$.  $m'$ is a monomial of $L_n(m'A_{n,b})$. By property (3), we can apply (3) of Lemma \ref{remont} and we get $M\in\mathcal{M}(L(m))$ dominant of the form $M=m'A_{n,b}A_{n-1,bq^{-1}}\cdots A_{n-r,bq^{-r}}$ with $r\geq 0$. From property (2), we have $M = m$. So $R=0$. So $n-r=K$, $bq^{-r}=a_Kq^{\lambda_K}$, that is to say $b=a_Kq^{\lambda_K+n-K}$, contradiction.\qed

\subsection{Type $B$}\label{secb} In this section \ref{secb}, we suppose that $\Glie$ is of type $B_n$.

\subsubsection{Preliminary results for type $B$}

\begin{lem}\label{etoileb} Let $a\in\CC^*$, $m\in\ZZ[Y_{i,aq^r}]_{i\in I,r\in\ZZ}$ a dominant monomial. Consider $m'\in\mathcal{M}(L(m))$ $\{1,\cdots,n-1\}$-dominant such that all $m''\in\mathcal{M}(L(m))$ satisfying $v(m''m^{-1}) < v(m'm^{-1})$ is thin. Suppose that $m'$ is not dominant and let $R=\text{min}\{r\in\ZZ|u_{n,aq^r}(m') < 0\}$. Then there is $M\in\mathcal{M}(L(m))$ $\{1,\cdots,n-1\}$-dominant such that $m\geq M > m'$, $m'M^{-1}\in\ZZ[A_{i,aq^{R+2(i-n)+4r-1}}]_{i\in I,r\leq 0}$, $u_{n,aq^R}(M) = 0$ and for all $r\leq R$, $u_{n,aq^r}(M)\geq 0$ and $\sum_{l\geq 0}u_{n,aq^{R-2-4l}}(M)>0$.\end{lem}

\demo For the shortness of notations, we suppose that $m'Y_{n,aq^R}$ is dominant (the proof is exactly the same if $R \neq \text{max}\{r\in\ZZ|u_{n,aq^r}(m') < 0\}$). First there is 
$$m_0=A_{n,aq^{R-1}}A_{n-1,aq^{R-3}}\cdots A_{n-\alpha,aq^{R-1-2\alpha}}m'\in\mathcal{M}(L(m)),$$
where $\alpha\geq 0$ and $m_0$ is $\{1,\cdots,n-1\}$-dominant. If $\alpha = 0$ we take $M = m_0$. Otherwise, $u_{n,aq^{R-4}}(m_0) = -1$ and $u_{n,b}(m_0) > 0$ implies $b=aq^{R-4}$. We continue and we get inductively (at each step the involved monomials are thin by assumption) :
$$m_r=A_{n,aq^{R-1-4r}}A_{n-1,aq^{R-3-4r}}\cdots A_{n-\alpha+r,aq^{R-1-2\alpha-2r}}m_{r-1}\in\mathcal{M}(L(m)),$$
where $1\leq r\leq \alpha$ and $m_r$ is $\{1,\cdots,n-1\}$-dominant. We take $M=m_{\alpha}$ and the properties are satisfied by construction.\qed

\begin{lem}\label{remontbii} Let $L(m)$ be a simple $\U_q(\Lo\Glie)$-module. Let $m'\in\mathcal{M}(L(m))$ such that all $m''\in\mathcal{M}(L(m))$ satisfying $v(m''m^{-1}) < v(m'm^{-1})$ is thin. Suppose that there are $j\in (I-\{n\})$ such that $u_{j,b}(m') < 0$ and $mY_{j,b}$ is dominant. Moreover we suppose that if $j\neq 1$, then $u_{j-1,bq^{-2}}(m') > 0$. Then there is $M\in\mathcal{M}(L(m))$ dominant satisfying one of the following conditions :

\begin{enumerate}
\item $M = m' A_{j,bq^{-2}}A_{j+1,bq^{-4}}\cdots A_{j+r,bq^{-2-2r}}$ where $0\leq r\leq n-j$,
 
\item $M = m'(A_{j,bq^{-2}}A_{j+1,bq^{-4}}\cdots A_{n-1,bq^{-2n+2j}})A_{n,bq^{-2n+2j}}M'$ where \\$M'\in\ZZ[A_{k,bq^{-2n+2j+2(k-n)-4l}}]_{k  < n,l\geq 0}\ZZ[A_{n,bq^{-2n+2j-4l}}]_{l\geq 1}$,

\item $M = m'(A_{j,bq^{-2}}A_{j+1,bq^{-4}}\cdots A_{n,bq^{-2-2n+2j}})A_{n,aq^{-2n+2j}}$,

\item $M = m'(A_{j,bq^{-2}}A_{j+1,bq^{-4}}\cdots A_{n,bq^{-2-2n+2j}})A_{n,bq^{-2n+2j}}A_{n-1,bq^{-2n-2+2j}}$.
\end{enumerate}

\noindent Moreover

in case (1), we have $u_{j+r,bq^{-2-2r-r_{j+r}}}(M)=1$,

in case (2), we have $u_{n-1,bq^{-2n+2j-2}}(M)=1$ and $\sum_{l\geq 0}u_{n,bq^{-2n+2j-1-4l}}(M) > 0$,

in case (3), we have $u_{n,bq^{-3-2n+2j}}(M)=u_{n,bq^{-1-2n+2j}}(M)=1$,

in case (4), we have $u_{n-1,bq^{-2n-4+2j}}(M)=1$.\end{lem}

\demo Thanks to the hypothesis $u_{j-1,bq^{-2}}(m') > 0$, we can suppose that $j=1$ . By using (3) of Lemma \ref{etoile} with $\Glie_{\{1,\cdots,n-1\}}$ of type $A_{n-1}$, we get 
$m_1 = m'A_{1,bq^{-2}}A_{2,bq^{-4}}\cdots A_{1+r,bq^{-2(r+1)}}\in\mathcal{M}(L(m)),$ 
$\{1,\cdots, n-1\}$-dominant.

If $m_1$ is dominant, then the condition (1) is satisfied, and we set $M=m_1$. 

Otherwise we have $r=n-2$, $u_{n-1,bq^{-2n}}(m_1)=1$, $m_1$ is not $n$-dominant and $(u_{n,d}(m_1)<0\Rightarrow d=bq^{-2n+3}\text{ or }d=bq^{-2n+1})$. 

If $u_{n,bq^{-2n+1}}(m_1)\geq 0$ and $u_{n,bq^{-2n+3}}(m_1) = -1$, then we can use Lemma \ref{etoileb} and so condition (2) is satisfied.

If $u_{n,bq^{-2n+1}}(m_1)=-1$ and $u_{n,bq^{-2n+3}}(m_1)\geq 0$, $M=m_1A_{n,bq^{-2n}}\in\mathcal{M}(L(m))$ is dominant, so condition (1) is satisfied. 

If $u_{n,bq^{-2n+1}}(m_1)=-1$ and $u_{n,bq^{-2n+3}}(m_1) = -1$, $m_2=m_1A_{n,bq^{-2n}}A_{n,bq^{-2n+2}}\in\mathcal{M}(L(m))$ is $n$-dominant and $(u_{n-1,d}(m_2)<0\Rightarrow d=bq^{-2n+2})$. If $m_2$ is dominant, condition (3) is satisfied. If $u_{n-1,bq^{-2n+2}}(m_2)= -1$, $M=m_2A_{n-1,bq^{-2n}}\in\mathcal{M}(L(m))$ is dominant as $u_{n,bq^{-2n+1}}(m_2)=u_{n,bq^{-2n-1}}(m_2)=1$. So condition (4) is satisfied.

The additional properties in the end of the statement are clear by construction of $M$.\qed

\subsubsection{Kirillov-Reshetikhin modules $W_{\lambda,a}^{(n)}$} Now we consider the case of a Kirillov-Reshetikhin module in the node $n$, that is to say a minimal affinization of $V(\lambda\Lambda_n)$ (observe that in this case condition (I) and condition (II) of Theorem \ref{class} are satisfied).

\begin{lem}\label{proofspebn} Let $m=X_{\lambda,a}^{(n)}$. Then
\begin{enumerate}
\item For all $m'\in \mathcal{M}(L(m))$ and $\mu\in\{\lambda,\lambda - 2,\cdots,-\lambda+2\}$, $v_{n,aq^\mu}(m'm^{-1})\geq 1$ implies $v_{n,aq^{\lambda}}(m'm^{-1})\geq 1,v_{n,aq^{\lambda -2}}(m'm^{-1})\geq 1,\cdots,v_{n,aq^{\mu}}(m'm^{-1})\geq 1$.

\item $L(m)$ is special.

\item $L(m)$ is thin.

\item Let $m'\in\mathcal{M}(L(m))$ satisfying $\sum_{r\in\ZZ,i < n}v_{i,aq^{\lambda+2n-2i+4r}}(m'm^{-1})\geq 0$. Let $j=\text{min}\{i<n|\sum_{r\in\ZZ}v_{i,aq^{\lambda+2n-2i+4r}}(m'm^{-1})\neq 0\}$. We have 
\begin{equation*}
\begin{split}
v_{j,aq^{\lambda+2n-2j}}(m'm^{-1})&=v_{j,aq^{\lambda+2n-2j-4}}(m'm^{-1})
\\&=\cdots=v_{j,aq^{\lambda+2n-2j-4R}}(m'm^{-1})=1,
\end{split}
\end{equation*}
where $R=\sum_{r\in\ZZ}v_{j,aq^{\lambda+2n-2j+4r}}(m'm^{-1})-1$.
\end{enumerate}\end{lem}

\demo (1) follows from Lemma \ref{geneun}. (2) follows from Theorem \ref{formerkr}. 

Let us prove (3). From property (2) and Proposition \ref{proofthin}, it suffices to prove that all monomials of $\mathcal{M}(L(m))$ are thin. From Lemma \ref{thinmon}, we can suppose that there is $m'\in\mathcal{M}(L(m))$ such that there are $l\in I, d\in\CC^*$ satisfying $u_{l,d}(m')=2$ and such that all $m''\in\mathcal{M}(L(m))$ satisfying $v(m''m^{-1}) < v(m'm^{-1})$ is thin. We distinguish three cases $(\alpha)$, $(\beta)$, $(\gamma)$.

$(\alpha)$ Suppose that there is $c\in\CC^*$ such that $u_{n,c}(m')\geq 2$. Then one of the two following condition is satisfied.

$(\alpha.i)$ : There is $b\in\CC^*$ such that $u_{n-1,b}(m') = -1$, $(u_{n-1,d}(m')<0\Rightarrow d=b)$, ($u_{n,bq^{-1}}(m')=2$ or $u_{n,bq^{-3}}(m')=2$) and $(u_{n,d}(m')=2\Rightarrow (d=bq^{-1}\text{ or }d=bq^{-3}))$.

$(\alpha.ii)$ : There is $b\in\CC^*$ such that $u_{n-1,b}(m') = u_{n-1,bq^2}(m') = -1$, $u_{n,bq^{-1}}(m')=2$, $(u_{n-1,d}(m')<0\Rightarrow (d = b\text{ or } d = bq^2))$ and $(u_{n,d}(m')=2\Rightarrow d=bq^{-1})$.

\noindent Otherwise, by using Proposition \ref{jdecomp}, we would get $m''\in\mathcal{M}(L(m))$ such that $v(m''m^{-1}) < v(m'm^{-1})$ and $m''$ does not satisfy property (3). 

First suppose that the condition $(\alpha.i)$ is satisfied. We have the following subcases :

$(\alpha.i.1)$ : $u_{n,bq^{-1}}(m')\geq 1$ and $u_{n,bq^{-3}}(m')\geq 1$. Then $m'A_{n-1,bq^{-2}}\in \mathcal{M}(L(m))$ is $(I-\{n-2\})$-dominant and by (3) of Lemma \ref{remont} with $\Glie_{\{1,\cdots,n-1\}}$ we get $M\in\mathcal{M}(L(m))$ dominant such that $u_{n-R,bq^{-2-2R}}(M)\geq 1$ for an $R\geq 1$. By property (2), $M = m$, contradiction.

$(\alpha.i.2)$ : $u_{n,bq^{-3}}(m')=2$ and $u_{n,bq^{-1}}(m')=0$. Then $m'' = m'A_{n-1,bq^{-2}}\in\mathcal{M}(L(m))$ and $Y_{n,bq^{-3}}Y_{n,bq^{-1}}^{-1}$ appears in $m''$. So by Lemma \ref{jdecomp} there is $m'''\in\mathcal{M}(L(m))$ such that $m'' < m'''$ and $u_{n,bq^{-3}}(m''')\geq 2$, contradiction. 

$(\alpha.i.3)$ : $u_{n,bq^{-1}}(m')=2$ and $u_{n,bq^{-3}}(m')=0$. Then 
$m''=m'A_{n-1,bq^{-2}}\in \mathcal{M}(L(m))$ is $(I-\{n-2,n\})$-dominant and $\forall j\in I, (u_{j,bq^{l}}(m'')>0\Rightarrow l\leq -2)$. So by Lemma \ref{plushaut} we get $M\in\mathcal{M}(L(m))$ dominant such that $v_{n,bq^{-4}}(m''M^{-1})\geq 1$ and $m''M^{-1}\in\ZZ[A_{j,bq^l}]_{j\in I, l\leq -3}$. So $u_{n,bq^{-1}}(M)=u_{n,bq^{-1}}(m'')=1$. By property (2), $M=m$. By property (1), we have $v_{n,b}(m''M^{-1})\geq 1$, contradiction.

Now we suppose that $(\alpha.ii)$ is satisfied. We have the following subcases :

$(\alpha.ii.1)$ : $u_{n,bq}(m') = 0$. Then $m''=m'A_{n-1,bq^{-2}}A_{n-1,b}A_{n,b}\in\mathcal{M}(L(m))$ and $Y_{n-1,bq^{-4}}Y_{n-1,bq^{-2}}Y_{n-1,b}^{-1}Y_{n,bq^{-1}}$ appears in $m''$. So $m'''=m''A_{n-1,bq^{-2}}\in\mathcal{M}(L(m))$ and $u_{n-1,bq^{-4}}(m''')=2$, contradiction.

$(\alpha.ii.2)$ : $u_{n,bq}(m') = u_{n,bq^{-3}}(m') = 1$. Consider $m''=m'A_{n-1,b}A_{n-1,bq^{-2}}$. Then $m''\in \mathcal{M}(L(m))$ is $(I-\{n-2\})$-dominant and by (3) of Lemma \ref{remont} with $\Glie_{\{1,\cdots, n-2\}}$ we get $M\in\mathcal{M}(L(m))$ dominant such that 
$$u_{n-r_1,bq^{-2r_1}}(M)\geq 1\text{ and }u_{n-r_2,bq^{-2-2r_2}}(M)\geq 1$$ 
with $r_1,r_2\geq 1$. By property (2), $M=m$, contradiction.

$(\alpha.ii.3)$ : $u_{n,bq}(m') = 1$ and $u_{n,bq^{-3}}(m') = 0$. Then $$m''=m'A_{n-1,bq^{-2}}A_{n-1,b}A_{n,bq^{-4}}\in\mathcal{M}(L(m))$$
is $(I-\{n-2\})$-dominant, $Y_{n-1,bq^{-2}}Y_{n,bq^{-5}}$ 
appears in $m''$ and $\forall j\in I, (u_{j,d}(m'')=-1\Rightarrow ((j,d) 
= (n-2,bq^{-2})\text{ or }(j,d)=(n-2,b)))$. So by (3) of Lemma \ref{remont}, there is $m'''\in\mathcal{M}(L(m))$
of the form $m'''=A_{n-2,bq^{-2}}A_{n-3,bq^{-4}}\cdots A_{n-R,bq^{2-2R}}m''$ with $R\geq 1$ such that $\forall j\in I, (u_{j,d}(m''')= -1\Rightarrow (j,d)=(n-2,bq^{-2}))$. We have $u_{n,bq^{-5}}(m''') = u_{n-R,bq^{-2R}}(m''') = 1$. If $m'''$ is dominant, we have $m'''=m$, contradiction. So $u_{n-2,bq^{-2}}(m''')=-1$. As moreover $u_{n,bq^{-5}}(m''')=1$, we have a dominant monomial $M\in\mathcal{M}(L(m))$ of the form :
\begin{equation*}
\begin{split}
M=&m'''(A_{n-2,bq^{-4}}A_{n-1,bq^{-6}}A_{n,bq^{-8}})(A_{n-3,bq^{-6}}A_{n-2,bq^{-8}}A_{n-1,bq^{-10}}A_{n,bq^{-12}})
\\&\cdots (A_{n-r,bq^{-2r}}A_{n-r+1,bq^{-2-2r}}\cdots A_{n,bq^{-4r}})
\\&\times (A_{n-r-1,bq^{-2-2r}}A_{n-r,bq^{-4-2r}}\cdots A_{n-r-1+r',bq^{-2-2r-2r'}}),
\end{split}
\end{equation*}
where $r\geq 1$ and $r+1\geq r'\geq 0$. By property (1), we have $M=m$. So we have $u_{n-R,bq^{-2R}}(m)=u_{n-R,bq^{-2R}}(m''')=1$, contradiction. 

$(\beta)$ Suppose that there is $b\in\CC^*$ such that $u_{n-1,b}(m')\geq 2$. Then we have $(u_{n-2,d}(m') < 0\Rightarrow d=bq^2)$ and $(u_{n,d}(m') < 0\Rightarrow d=bq)$. By (3) of Lemma \ref{remont} with $J=\{1,\cdots, n-1\}$ and $J=\{n\}$, we get a dominant monomial $M\in\mathcal{M}(L(m))$ satisfying one of the two following condition :

$(\beta.1)$ $u_{j_1,bq^{2j_1-2n+2}}(M) = 1$, $u_{n,bq^{-1}}(M) = 1$ with $j_1 \leq n-1$. 

$(\beta.2)$ $u_{j_1,bq^{2j_1-2n+2}}(M) = 1$, $u_{n-1,b}(M) = 1$ with $j_1 \leq n-2$.

\noindent From (2) we have $m=M$, contradiction.

$(\gamma)$ Suppose that there is $i \leq n-2$ and $b\in\CC^*$ such that $u_{i,b}(m')\geq 2$. Then $m'$ is $(\{1,\cdots,i-2\}\cup \{i\}\cup \{i+2,\cdots,n\})$-dominant. We have $(u_{i-1,d}(m')<0\Rightarrow d=bq^2)$ and  $(u_{i+1,d}(m')<0\Rightarrow d=bq^2)$. By applying (3) of Lemma \ref{remont} and Lemma \ref{remontbii}, we get $M\in\mathcal{M}(L(m))$ dominant such that $(M)^{\{1,\cdots,n-2\}}\neq 1$. From property (2) we have $m=M$, contradiction.

Now we prove property (4) by induction on $v(m'm^{-1})\geq 0$. Let $j$ be as in property (4). For $v(m'm^{-1})=0$ we have $m'=m$ and the result is clear. We suppose that property (4) is satisfied for $m''$ satisfying $v(m''m^{-1}) < v(m'm^{-1})$. Let $R\geq 0$ maximal such that $m'm^{-1}\leq A_{j,aq^{\lambda+2n-2j}}^{-1}A_{j,aq^{\lambda+2n-2j-4}}^{-1}\cdots A_{j,aq^{\lambda+2n-2j+4-4R}}^{-1}$. We suppose moreover that 
$$m'm^{-1}\leq A_{j,aq^{\lambda+2n-2j}}^{-1}A_{j,aq^{\lambda+2n-2j-4}}^{-1}\cdots A_{j,aq^{\lambda+2n-2j+4-4R}}^{-1}A_{j,b}^{-1},$$ 
with $b = aq^{\lambda+2n-2j-4\mu}$ and $\mu\neq R$. If $j\geq 2$, we have $u_{j-1,b}(m')=1$. By the induction hypothesis on $v$, $m'$ is $(I-\{j\})$-dominant, $u_{j,bq^2}(m')=-1$ and ($u_{j,d}(m')<0\Rightarrow d = bq^2$). By property (3), we can apply Lemma \ref{remontbii} and we get a dominant monomial $M\in\mathcal{M}(L(m))$. From property (2), we have $M = m$. As $u_{n-1}(m)=0$, we are in the situation (1) or (3) of Lemma \ref{remontbii}. So
$$m(m')^{-1}=A_{j,b}A_{j+1,bq^{-2}}\cdots A_{n-1,bq^{-2(n-1-j)}}A_{n,bq^{-2(n-j)}}(A_{n,bq^{2-2(n-j)}})^{\epsilon},$$ 
where $\epsilon\in\{0,1\}$. So $b=aq^{\lambda+2n-2j}$, $\mu=0$ and $R=0$, contradiction.\qed

\subsubsection{Condition (I)} Now we treat the general case of minimal affinization satisfying condition (I) of Theorem \ref{class} (except the Kirillov-Reshetikhin modules $W_{k,a}^{(n)}$ already studied in Lemma \ref{proofspebn}).

\begin{lem}\label{proofspeb} Let $\lambda\in P^+$ and $L(m)$ be a minimal affinization of $V(\lambda)$ such that $m$ satisfies condition (I) of Theorem \ref{class}. Let $K=\text{min}\{i\in I|\lambda_i\neq 0\}$. We suppose that $K\leq n-1$. Then the following conditions are satisfied :
\begin{enumerate}
\item For all $m'\in \mathcal{M}(L(m))$ satisfying $\sum_{r\in\ZZ}v_{K,a_Kq^{2\lambda_K+4r}}(m'm^{-1})\geq 1$, we have $v_{K,a_Kq^{2\lambda_K}}(m'm^{-1})\geq 1$.

\item $L(m)$ is special.

\item $L(m)$ is thin.

\item Let $m'\in\mathcal{M}(L(m))$ satisfying $\sum_{r\in\ZZ,i<n}v_{i,a_iq^{2\lambda_i+4r}}(m'm^{-1})\geq 1$. We have 
\begin{equation*}
\begin{split}
v_{j,a_kq^{r_k\lambda_k+2k-2j}}(m'm^{-1})&=v_{j,a_kq^{r_k\lambda_k+2k-2j-4}}(m'm^{-1})
\\&= \cdots =v_{j,a_kq^{r_k\lambda_k+2k-2j-4R}}(m'm^{-1})=1,
\end{split}
\end{equation*}
where 
$$j=\text{min}\{i<n|\sum_{r\in\ZZ}v_{i,a_iq^{2i+4r}}(m'm^{-1})\geq 1\},$$ 
$k=\text{min}\{i\geq j|\lambda_i\neq 0\}$ and $R=(\sum_{r\in\ZZ}v_{j,a_jq^{2\lambda_j+4r}}(m'm^{-1}))-1$.

\item Let $m'\in\mathcal{M}(L(m))$ such that $\text{min}\{i|\sum_{r\in\ZZ}v_{i,a_iq^{2\lambda_i+4r}}(m'm^{-1})\geq 1\}=n$. Then $\lambda_n\neq 0$ and  
\begin{equation*}
\begin{split}
v_{n,a_nq^{\lambda_n}}(m'm^{-1})&=v_{n,a_nq^{\lambda_n-2}}(m'm^{-1})
\\&=\cdots=v_{n,a_nq^{\lambda_n-2R}}(m'm^{-1})=1,
\end{split}
\end{equation*}
where $R=\sum_{r\in\ZZ}v_{n,a_nq^{2\lambda_n+2r}}(m'm^{-1})-1$.
\end{enumerate}
\end{lem}

\demo We prove by induction on $u(m)\geq 0$ simultaneously that (1), (2), (3), (4) and (5) are satisfied.

For $u(m)=0$ the result is clear. Suppose that $u(m)\geq 1$. 

First we prove (1) by induction on $v(m'm^{-1})\geq 0$. For $v(m'm^{-1})=0$ we have $m'=m$ and the result is clear. In general suppose that for $m''$ such that $v(m''m^{-1}) < v(m'm^{-1})$ the property is satisfied, that $v_{K,a_Kq^{2\lambda_K}}(m'm^{-1}) = 0$ and $\sum_{r\in\ZZ}v_{K,a_Kq^{2\lambda_K+4r}}(m'm^{-1})\geq 1$. Observe that it follows from Lemma \ref{jdecomp} and corollary \ref{jdecompprime} that $m'$ is $(I-\{K\})$-dominant and $(m')^{\rightarrow (K,a_Kq^{2\lambda_K})}$ is dominant. 

If $m'$ is not dominant, by corollary \ref{jdecompprime}, there is $m''\in\mathcal{M}(L(m))$ $K$-dominant such that $m'$ is a monomial of $L_{K,a_Kq^{2\lambda_K-2}}(m'')$. Moreover from Proposition \ref{aidesldeux}, there is $b\in a_Kq^{2\lambda_K+4\ZZ}$ such that $A_{K,b}m'\in\mathcal{M}(L(m))$. By the induction property on $v$, we have $\sum_{r\in\ZZ}v_{K,a_Kq^{2\lambda_K+4r}}(m'A_{K,b}m^{-1})=0$. So $m''=m'A_{K,b}$. But $m''\in m^{\rightarrow (K)}\mathcal{M}(L(m(m^{\rightarrow (K)})^{-1}))$. As $u(m(m^{\rightarrow (K)})^{-1}) < u(m)$, we have property (4) for $L(m(m^{\rightarrow (K)})^{-1})$ and we get 
$$(m'')^{\rightarrow (K)} \in Y_{K,a_Kq^{2\lambda_K-2}}Y_{K,a_Kq^{2\lambda_K-6}}\cdots Y_{K,a_Kq^{-2\lambda_K+2-4R'}}\ZZ[Y_{K,a_Kq^{4r+2\lambda_K}}]_{r\in\ZZ}$$
with $R'\geq 0$. By Lemma \ref{aidesldeux}, $m'$ is not a monomial of $\mathcal{M}(L_K(m''))$, contradiction. So $m'$ is dominant.

Let us prove that $\sum_{r\in\ZZ}v_{K,a_Kq^{2\lambda_K+4r+2}}(m'm^{-1}) = 0$. Observe that 
$$m'Y_{K,a_Kq^{2\lambda_K - 2}}^{-1}\in \mathcal{M}(L(mY_{K,a_Kq^{2\lambda_K-2}}^{-1})).$$ 
Moreover $u_{j,a}(m'(m^{\rightarrow (K)})^{-1}) < 0$ implies $j=K$ and $a = a_Kq^{2\lambda_K - 2}$. As we have $u(mY_{K,a_Kq^{2\lambda_K-2}}^{-1}) < u(m)$, properties (2) and (3) are satisfied by $L(mY_{K,a_Kq^{2\lambda_K-2}}^{-1})$. So we can use (2) of Lemma \ref{stara} for $\Glie_{\{1,\cdots,n-1\}}$ of type $A_{n-1}$ and we get a monomial
$$m''\in\mathcal{M}(L(mY_{K,a_Kq^{2\lambda_K-2}}^{-1}))\cap (m'Y_{K,a_Kq^{2\lambda_K-2}}^{-1}\ZZ[A_{j,a_Kq^{2\lambda_K+4r+2(K-j)}}]_{j\leq n-1,r\in\ZZ}),$$ 
which is $\{1,\cdots,n-1\}$-dominant and satisfying $v_{n-1}(m'Y_{K,a_Kq^{2\lambda_K-2}}^{-1}(m'')^{-1})\leq 1$. If $v_{n-1}(m'Y_{K,a_Kq^{2\lambda_K-2}}^{-1}(m'')^{-1}) = 0$, then $mY_{K,a_Kq^{2\lambda_K-2}}^{-1}=m''$ and the result is clear. Otherwise, consider the unique $b\in a_Kq^{2\lambda_K+2(n-K)+4\ZZ}$ such that $v_{n-1,b}(m'Y_{K,a_Kq^{2\lambda_K-2}}^{-1}(m'')^{-1}) = 1$. We have $u_{n,d}(m'') < 0\Rightarrow (d = bq\text{ or } d = bq^{-1})$. If $u_{n,bq}(m'')= 0$ we use Lemma \ref{etoileb} and we get the result. If $u_{n,bq}(m'') = -1$ and $u_{n,bq^{-1}}(m'') = 0$, we use Lemma \ref{etoileb}, and in particular we get a monomial 
$$(m Y_{K,a_Kq^{2\lambda_K-2}}^{-1}) A_{j,d}^{-1}\in \mathcal{M}(L(mY_{K,a_Kq^{2\lambda_K-2}}^{-1}))$$ 
where $d\notin a_Kq^{2\lambda_K+4\ZZ}$, contradiction with condition (II). In the same way if we have $u_{n,bq}(m'') = -1$ and $u_{n,bq^{-1}}(m'') = -1$, then we get a contradiction by using twice Lemma \ref{etoileb}.

Now it suffices to prove that the conditions of Proposition \ref{racourc} with $i = K$ are satisfied. 

Condition (i) of Proposition \ref{racourc} : if $M>m'$ is in $\mathcal{M}(L(m))$, we have necessarily $v_{K,a_Kq^{2\lambda_K}}(Mm^{-1}) = 0$. So by induction hypothesis $\sum_{r\in\ZZ}v_{K,a_Kq^{2\lambda_K+4r}}(Mm^{-1}) = 0$, and so $v_K(Mm^{-1})=0$. So if we suppose moreover that $M\in m'\ZZ[A_{K,a}]_{a\in\CC^*}$, we have necessarily $M=m\prod_{a\in\CC^*}A_{K,a}^{v_{K,a}(m'm^{-1})}$, and so we get the uniqueness. For the existence, it suffices to prove that this $M=m\prod_{a\in\CC^*}A_{K,a}^{v_{K,a}(m'm^{-1})}$ is in $\mathcal{M}(L(m))$. By Lemma \ref{oudeux}, there is $j\in I$, $M'\in\mathcal{M}(L(m))$ $j$-dominant such that $M'>m'$ and $M'\in m'\ZZ[A_{j,a}]_{a\in\CC^*}$. By induction hypothesis on $v$ we have $j=K$, and so by uniqueness $M'=M$. 

Condition (ii) of Proposition \ref{racourc} : by construction of $M$ we have $v_K(Mm^{-1})=0$. 

Condition (iii) of Proposition \ref{racourc} : first observe that 
$$M\in m^{\rightarrow (K)}\mathcal{M}(L(m(m^{\rightarrow (K)})^{-1})).$$ 
As $u(m(m^{\rightarrow (K)})^{-1}) < u(m)$, we have property (4) for $L(m(m^{\rightarrow (K)})^{-1})$ and we get 
$$(M)^{\rightarrow (K)} \in Y_{K,a_Kq^{2\lambda_K-2}}Y_{K,a_Kq^{2\lambda_K-6}}\cdots Y_{K,a_Kq^{-2\lambda_K+2-4R'}}\ZZ[Y_{K,a_Kq^{4r+2\lambda_K}}]_{r\in\ZZ}$$
with $R'\geq 0$. By Lemma \ref{aidesldeux}, $m'$ is not a monomial of $\mathcal{M}(L_K(M))$. 

Condition (iv) of Proposition \ref{racourc} : consider a monomial $m''\in\mathcal{M}(\U_q(\Lo\Glie_K).V_M)$ such that $v(m''m^{-1})<v(m'm^{-1})$. We have $m''\leq MA_{K,a_Kq^{2\lambda_k}}^{-1}\ZZ[A_{K,d}^{-1}]_{d\in\CC^*}$ and so $(m'')^{\rightarrow (K,a_Kq^{2\lambda_K})}$ is right negative, so $m''$ is not $K$-dominant. 

Condition (v) of Proposition \ref{racourc} : clear by the induction property on $v$.

Now we prove (2). Let $J=\{i\in I|K<i\}$. From Lemma \ref{produit}, 
$$\mathcal{M}(L(m))\subset (m^{\rightarrow (J)}\mathcal{M}(L(m^{\rightarrow (K)})))\cup(\mathcal{M}(L(m^{\rightarrow (J)}))m^{\rightarrow (K)}).$$ 
As all monomials of $m^{\rightarrow (J)}(\chi_q(L(m^{\rightarrow (K)}))-m^{\rightarrow (K)})$ are lower than $mA_{K,a_Kq^{\lambda_K}}^{-1}$ (Theorem \ref{formerkr}) which is right-negative, they are not dominant. Let $m'$ be a monomial in $(\mathcal{M}(L(m^{\rightarrow (J)}))m^{\rightarrow (K)}-\{m\})$. As $u(m^{\rightarrow (J)}) < u(m)$, the induction property implies that $m'(m^{\rightarrow (K)})^{-1}$ is not dominant. If $\sum_{l\in\ZZ}v_{K,a_Kq^{2\lambda_K+4l}}(m'm^{-1})\geq 1$, it follows from property (1) that $m'$ is lower than $mA_{K,a_Kq^{2\lambda_K}}^{-1}$ which is right-negative, so $m'$ is not dominant. We suppose that $\sum_{l\in\ZZ}v_{K,a_Kq^{2\lambda_K+4l}}(m'm^{-1})=0$. We have for all $l\in\ZZ$, $u_{K,q^{2K+4l}}(m'(m^{\rightarrow (K)})^{-1})\geq 0$, and so there is $(i,a)\in I\times \CC^*$ not of the form $(K,q^{2K+4l})$ with $l\in\ZZ$ such that $u_{i,a}(m'(m^{\rightarrow (K)})^{-1})<0$. So $$u_{i,a}(m')=u_{i,a}(m'(m^{\rightarrow (K)})^{-1})<0$$ 
and $m'$ is not dominant. So $L(m)$ is special.

Now we prove (3). From property (2) and Proposition \ref{proofthin}, it suffices to prove that all monomials of $\mathcal{M}(L(m))$ are thin. From Lemma \ref{thinmon}, we can suppose that there is $m'\in\mathcal{M}(L(m))$ such that there are $l\in I, d\in\CC^*$ satisfying $u_{l,d}(m')=2$ and such that all $m''\in\mathcal{M}(L(m))$ satisfying $v(m''m^{-1}) < v(m'm^{-1})$ is thin. We consider subcases as in the proof of Lemma \ref{proofspebn}. 

If $(\alpha.i.1)$ is satisfied, we get $u_{n-R,bq^{-2-2R}}(m)\geq 1$ with $R\geq 1$ and ($u_{n,bq^{-1}}(m)\geq 1$ or $u_{n,bq^{-3}}(m)\geq 1$). As $-2-2R < -3$, we get a contradiction with condition (I) of Theorem \ref{class}.

If $(\alpha.i.2)$ is satisfied, we get a contradiction as for Lemma \ref{proofspebn}.

If $(\alpha.i.3)$ is satisfied, as for Lemma \ref{plusgrand} we get $m''\in\mathcal{M}(L(m))\cap m\ZZ[A_{i,bq^l}^{-1}]_{i\in I, l\leq -3}$ such that $u_{n,bq^{-1}}(m)=u_{n,bq^{-1}}(m'')=1$, $v_{n,bq^{-4}}(m''m^{-1})\geq 1$. From Lemma \ref{plusgrand} and Lemma \ref{produit}, we have $m''\in m^{\{1,\cdots, n-1\}}\mathcal{M}(L(m^{\rightarrow (n)}))$, and we get a contradiction as for Lemma \ref{proofspebn}.

If condition $(\alpha.ii.1)$ is satisfied, we get a contradiction as for Lemma \ref{proofspebn}.

If condition $(\alpha.ii.2)$ is satisfied, we get as in the proof of Lemma \ref{proofspebn} that $$u_{n-r_1,bq^{-2r_1}}(m)\geq 1\text{ and }u_{n-r_2,bq^{-2-2r_2}}(m)\geq 1$$ 
with $r_1,r_2\geq 1$. Contradiction with condition (I) of Theorem \ref{class}.

If condition $(\alpha.ii.3)$ is satisfied : we follow the proof of Lemma \ref{proofspebn} and we get $m'''$. If $m'''$ is dominant, we have $u_{n,bq^{-5}}(m''')=u_{n-R,bq^{-2R}}(m''')=1$ with $-2R-(-5)\leq 3 < 2(n-(n-R))+4$, contradiction with condition (I) of Theorem \ref{class}. So $m'''$ is not dominant. Let $R,r,r'$ and $M$ dominant defined in the proof of Lemma \ref{proofspebn}. From the property (2) we have $m=M$. Observe that $r'\leq r+1$. We have $u_{n-R,bq^{-2R}}(m)=1$. We study two cases :

if $n-r-1+r'=n$, we have moreover $u_{n,bq^{-3-2r-2r'}}(m) = 1$. But $(-3-2r-2r')-(-2R)\leq 2R-4 < 2(n-(n-R))$, contradiction with condition (I) of Theorem \ref{class}.

if $n-r-1+r'\leq n-1$, we have moreover $u_{n-r-1+r',bq^{-4-2r-2r'}}(m) = 1$. Let $d = (n-r-1+r'-(n-R)) = - 1 + (R+r'-r)$ and $D = (-4-2r-2r') + 2R = 2d   - 4r' - 2$. 
\\If $d < 0$, condition (I) implies $D\geq -2d + 4$, so $0\leq D+2d-4 = 4d-4r'-6  < 0$, contradiction. 
\\If $d = 0$, condition (I) implies $D \in 4\ZZ$, contradiction as $D=-4r' - 2$. 
\\If $d > 0$, condition (I) implies $D \leq -4 - 2d$, so $0\geq D + 4 + 2d = 4d + 2 - 4r'= -2 + 4R - 4r$  and $n-r-1 < n-R < n-R-1+r'$. So the product $A_{n-r-1+r',bq^{-2-2r-2r'}}\cdots A_{n-r-1,bq^{-2-2r}}$ can not appear in $m(m''')^{-1}$ (for example we may use Theorem \ref{racourc} as in the proof of Lemma \ref{proofspe}), contradiction.

Now we suppose that there is $b\in\CC^*$ such that $u_{n-1,b}(m')\geq 2$. By property (2), we get as in the proof of Lemma \ref{proofspebn} that $m$ satisfies property $(\beta.1)$ or $(\beta.2)$ of Lemma \ref{proofspebn}. For $(\beta.1)$, we have $(2j_1-2n+2-(-1))= 2(j_1 - n)+3 < 2(n-j_1)+5$, contradiction with condition (I) of Theorem \ref{class}. For $(\beta.2)$, we have $(2j_1-2n+2-0)= 2(j_1-n)+2 < 2(n-1-j_1)+6$, contradiction with condition (I) of Theorem \ref{class}. 

Finally we suppose that there are $i\leq n-2$, $b\in\CC^*$ such that $u_{i,b}(m')\geq 2$. Then $m'$ is $(\{1,\cdots,i-2\}\cup \{i\}\cup \{i+2,\cdots,n\})$-dominant. We have $(u_{i-1,d}(m')<0\Rightarrow d=bq^2)$, and  $(u_{i+1,d}(m')<0\Rightarrow d=bq^2)$. By applying (3) of Lemma \ref{remont} and Lemma \ref{remontbii} (with $bq^2$ instead of $b$ and $i+1$ instead of $j$), we get a dominant monomial $M\in\mathcal{M}(L(m))$ satisfying one of the conditions 

$(\gamma.1)$ (case (1) of Lemma \ref{remontbii}) : $u_{j_1,bq^{2j_1-2i}}(M)\geq 1$, $u_{j_2,bq^{2i-2j_2+2-r_{j_2}}}(M)\geq 1$ with $j_1 < j_2$, $j_1 \leq i\leq j_2\leq n$, 

$(\gamma.2)$ (case (2) of Lemma \ref{remontbii}) : $u_{j_1,bq^{2j_1-2i}}(M)\geq 1$ and $u_{n-1,bq^{-2n+2i+2}}(M)\geq 1$ with $j_1 \leq i$,

$(\gamma.3)$ (case (3) of Lemma \ref{remontbii}) : $u_{j_1,bq^{2j_1-2i}}(M)\geq 1$ and $u_{n,bq^{1-2n+2i}}(M)\geq 1$, $u_{n,bq^{3-2n+2i}}(M)\geq 1$ with $j_1 \leq i$,

$(\gamma.4)$ (case (4) of Lemma \ref{remontbii}) : $u_{j_1,bq^{2j_1-2i}}(M)\geq 1$ and $u_{n-1,bq^{-2n+2i}}(M)\geq 1$ with $j_1\leq i$. 

\noindent From property (2) we have $m=M$. For $(\gamma.1)$, we have $2j_1-2i-(2i-2j_2+2-r_{j_2}) \leq 2(j_1+j_2) -4i \leq 2(j_2-j_1)$, contradiction with condition (I) of Theorem \ref{class}. For $(\gamma.2)$, we have $2j_1-2i-(-2n+2i+2)\leq 2(n-1-j_1) $, contradiction with condition (I) of Theorem \ref{class}. For $(\gamma.3)$, we have $2j_1-2i-(3-2n+2i)\leq 2(n-j_1)$, contradiction with condition (I) of Theorem \ref{class}. For $(\gamma.4)$, we have $2j_1-2i-(2i-2n) < 2(n-1-j_1)+4$, contradiction with condition (I) of Theorem \ref{class}.

Now we prove property (4) by induction on $v(m'm^{-1})\geq 0$. Let $j$ be as in property (4). For $v(m'm^{-1})=0$ we have $m'=m$ and the result is clear. We suppose that property (4) is satisfied for $m''$ such that $v(m''m^{-1}) < v(m'm^{-1})$. Let $R\geq 0$ maximal such that 
$$m'm^{-1}\leq A_{j,a_kq^{r_k\lambda_k+2(k-j)}}^{-1}A_{j,a_kq^{r_k\lambda_k+2(k-j)-4}}^{-1}\cdots A_{j,a_kq^{r_k\lambda_k+2(k-j)+4-4R}}^{-1}.$$ 
We suppose moreover that 
$$m'm^{-1}\leq A_{j,a_kq^{r_k\lambda_k+2(k-j)}}^{-1}A_{j,a_kq^{r_k\lambda_k+2(k-j)-4}}^{-1}\cdots A_{j,a_kq^{r_k\lambda_k+2(k-j)+4-4R}}^{-1}A_{j,b}^{-1}$$ 
with $b = a_kq^{r_k\lambda_k+2(k-j)-4\mu}$ and $\mu\neq R$. By the induction hypothesis on $v$, $m'$ is $(I-\{j\})$-dominant, $u_{j,bq^2}(m')=-1$ and ($u_{j,d}(m')<0\Rightarrow d = bq^2$). By property (3), we can apply Lemma \ref{remontbii} and we get a dominant monomial $M\in\mathcal{M}(L(m))$. From property (2), we have $M = m$. So we have one of the following situations :

Case (1) of Lemma \ref{remontbii} : $m = m' A_{j,b}A_{j+1,bq^{-2}}\cdots A_{j+r,bq^{-2r}}$ where $0\leq r\leq n-j$, and $u_{j+r,bq^{-2r-r_{j+r}}}(M)=1$. So $R=0$, $j+r=k$, $b=a_kq^{r_k\lambda_k+2r}=a_kq^{r_k\lambda_k+2(k-j)}$, contradiction.
 
Case (2) of Lemma \ref{remontbii} : $m = m'(A_{j,b}A_{j+1,bq^{-2}}\cdots A_{n-1,bq^{2-2n+2j}})A_{n,bq^{2-2n+2j}}M'$ where 
$$M'\in\ZZ[A_{p,bq^{2-2n+2j+2(p-n)-4l}}]_{p  < n,l\geq 0}\ZZ[A_{n,bq^{2-2n+2j-4l}}]_{l\geq 1},$$ 
and $u_{n-1,bq^{-2n+2j}}(m)=1$, so $b\in a_{n-1}q^{2\lambda_{n-1}+2(j-n)+2+4\ZZ}$. There is $l\geq 0$ such that $bq^{-2n+2j+1-4l} = a_nq^{\lambda_n - 1}$. So $b\in a_nq^{\lambda_n-2+2n-2j + 4\ZZ} = a_{n-1}q^{2\lambda_{n-1}+2(j-n)+4\ZZ}$ from condition (I) of Theorem \ref{class}, contradiction.

Case (3) of Lemma \ref{remontbii} : $u_{n,bq^{-1-2n+2j}}(m)=u_{n,bq^{1-2n+2j}}(m)=1$ and
$$m=m'(A_{j,bq}A_{j+1,bq^{-2}}\cdots A_{n,bq^{-2n+2j}})A_{n,bq^{2-2n+2j}}.$$ 
So $R=0$ and $b=a_nq^{\lambda_n+2n-2j-2}$. From condition (I) of Theorem \ref{class}, $a_nq^{\lambda_n}=a_kq^{r_k\lambda_k+2(k-n)+4r}$ with $r\in\ZZ$. So $b=a_kq^{r_k\lambda_k+2(k-j)+4r-2}$ is not of the form $a_kq^{r_k\lambda_k+2k-2j-4\mu}$, contradiction.

Case (4) of Lemma \ref{remontbii} : $u_{n-1,bq^{-2n-2+2j}}(m)=1$ and 
$$m=m'(A_{j,b}A_{j+1,bq^{-2}}\cdots A_{n,bq^{-2n+2j}})A_{n,bq^{2-2n+2j}}A_{n-1,bq^{-2n+2j}}.$$ 
So $bq^{-2n+2j}=a_{n-1}q^{2\lambda_{n - 1}}$, and so $b=a_{n-1}q^{2\lambda_{n-1}+2(n-j)}\in a_kq^{2\lambda_k+2(k-j)+2+4\ZZ}$, contradiction.

Now we prove property (5) by induction on $v(m'm^{-1})\geq 0$. For $v(m'm^{-1})=0$ we have $m'=m$ and the result is clear. We suppose that property (5) is satisfied for $m''$ such that $v(m''m^{-1}) < v(m'm^{-1})$ and we suppose that 
$$\text{min}\{i|\sum_{r\in\ZZ}v_{i,a_iq^{2\lambda_i+4r}}(m'm^{-1})\geq 1\}=n.$$ 
Let $R\geq 0$ maximal such that 
$$m'm^{-1}\leq A_{n,a_nq^{\lambda_n}}^{-1}A_{n,a_nq^{\lambda_n-2}}^{-1}\cdots A_{n,a_nq^{\lambda_n+2-2R}}^{-1}.$$ 
We suppose moreover that 
$$m'm^{-1}\leq A_{n,a_nq^{\lambda_n}}^{-1}A_{n,a_nq^{\lambda_n-2}}^{-1}\cdots A_{n,a_nq^{\lambda_n+2-2R}}^{-1}A_{n,b}^{-1}$$ 
with $b = a_nq^{\lambda_n-2\mu}$ and $\mu\neq R$. By the induction hypothesis on $v$, $m'$ is $(I-\{n\})$-dominant, $u_{n,bq}(m')=-1$ and ($u_{n,d}(m')<0\Rightarrow d = bq$). So $m''=A_{n,b}m'\in\mathcal{M}(L(m))$ is $(I-\{n-1\})$-dominant and $(u_{n-1,d}(m'')<0\Rightarrow d=b)$. If $u_{n-1,b}(m'')\geq 0$, $m''$ is dominant equal to $m$, so $R=0$ and $b=a_nq^{\lambda_n}$, contradiction. So $u_{n-1,b}(m'')<0$, $m''A_{n-1,bq^{-2}}\in\mathcal{M}(L(m))$ and $v_{n-1,bq^{-2}}(m'm^{-1})\geq 1$. So $bq^{-2}\notin a_{n-1}q^{2\lambda_{n-1}+4\ZZ}$ and $b\notin a_nq^{\lambda_n+4\ZZ}$. By lemma \ref{etoileb} there is $l\in\ZZ$ such that $bq^{-1-4l} = a_nq^{\lambda_n - 1}$, so $b\in a_nq^{\lambda_n+4\ZZ}$, contradiction.\qed

\subsubsection{Condition (II)} We study the general case of study condition (II) of Theorem \ref{class}. 

\begin{lem}\label{proofspebdeux} Let $\lambda\in P^+$ and $L(m)$ be a minimal affinization of $V(\lambda)$ such that $m$ satisfies condition (II) of Theorem \ref{class}. Let $K=\text{max}\{i\in I|\lambda_i\neq 0\}$. Then

(1) For all $m'\in\mathcal{M}(L(m))$, if $v_K(m'm^{-1})\geq 1$, then $v_{K,a_Kq_K^{\lambda_K}}(m'm^{-1})\geq 1$.

(2) $L(m)$ is special.

(3) $L(m)$ is thin.

(4) For all $m'\in \mathcal{M}(L(m))$ such that $v_n(m'm^{-1})=0$ we have 
\begin{equation*}
\begin{split}
v_{j,a_kq^{2(\lambda_k+j-k)}}(m'm^{-1}) & = v_{j,a_kq^{2(\lambda_k+j-k-2)}}(m'm^{-1})
\\&=\cdots=v_{j,a_kq^{2(\lambda_k+j-k-2R)}}(m'm^{-1})=1,
\end{split}
\end{equation*}
where $j=\text{max}\{i|v_i(m'm^{-1})\neq 0\}$, $k=\text{max}\{i\leq j|\lambda_i\neq 0\}$ and $R=v_j(m'm^{-1})-1$.
\end{lem}

Observe that Lemma \ref{proofspebn}, Lemma \ref{proofspeb} and Lemma \ref{proofspebdeux} combined with corollary \ref{simplification} imply Theorem \ref{fora} and Theorem \ref{egalun} for type $B$.

In this case we do not need to prove simultaneously the different properties.

\demo Property (4) : as $v_n(m'm^{\rightarrow (-1)})= 0$, it follows from Lemma \ref{depart} that $m'$ appears in $L_{\{1,\cdots n-1\}}(m)$. As $\Glie_J$ is of type $A_{n-1}$, the result is exactly property (4) of Lemma \ref{proofspe}.

Property (1) and (2) : as property (4) is satisfied, we can use the proof of property (1) and (2) of Lemma \ref{proofspe}.

Property (3) : the monomial $M=\prod_{i\in I, a\in\CC^*}Y_{i,a^{-1}q^{-r^{\vee}h^{\vee}}}^{u_{i,a}(m)}$ satisfies (I), and so it follows from Lemma \ref{proofspebn} and Lemma \ref{proofspeb} that $L(M)$ is thin. But from corollary \ref{invmon}, $\sigma^*L(M)\simeq L(m)$, and so we have property (3) (Lemma \ref{involchiq}).\qed

\subsection{Type $G_2$}

In this section we suppose that $\Glie$ is of type $G_2$. 

\begin{lem}\label{proofspeg} Let $m$ be a dominant monomial satisfying condition (I) of Theorem \ref{class}. Then $L(m)$ is special.\end{lem}

\demo From Lemma \ref{produit}, $\mathcal{M}(L(m))\subset \mathcal{M}(L(m^{\rightarrow (1)})\mathcal{M}(L(m^{\rightarrow (2)}))$. From Lemma \ref{formerkr}, if $m'$ is in $(\mathcal{M}(L(m^{\rightarrow (1)})-\{m^{\rightarrow (1)}\})\mathcal{M}(L(m^{\rightarrow (2)}))$, then $m'\leq mA_{1,a_1q^{3\lambda_1}}^{-1}$ which is right-negative, and so $m'$ is not dominant. Consider $m'=m^{\rightarrow (1)}m_2'$ where $m_2'\in(\mathcal{M}(L(m^{\rightarrow (2)}))-\{m^{\rightarrow (2)}\})$. It follows from Theorem \ref{formerkr} that $m_2'$ is right-negative. Suppose that $m'$ is dominant. In particular $m_2'$ is $2$-dominant and ($u_{1,b}(m_2')<0\Rightarrow (u_{1,b}(m_2')=-1\text{ and }b\in\{a_1q^{3-3\lambda_1},a_1q^{9-3\lambda_1},\cdots,a_1q^{3\lambda_1-3}\})$). From Lemma \ref{produit}, $m_2'\in \mathcal{M}(V_2(a_2q^{1-\lambda_2}))\mathcal{M}(V_2(a_2q^{3-\lambda_2}))\cdots \mathcal{M}(V_2(a_2q^{\lambda_2-1}))$. But for $b\in\CC^*$, it follows from \cite[Section 8.4.1]{her02} (with $1$ instead of $2$ and $2$ instead of $1$) that 
\begin{equation*}
\begin{split}
\chi_q(V_2(b)) = 
&Y_{2,b}+Y_{2,bq^2}^{-1}Y_{1,bq}+Y_{1,bq^7}^{-1}Y_{2,bq^4}Y_{2,bq^6}+Y_{2,bq^4}Y_{2,bq^8}^{-1}
+Y_{1,bq^5}Y_{2,bq^6}^{-1}Y_{2,bq^8}^{-1}
\\&+Y_{1,bq^{11}}^{-1}Y_{2,bq^{10}}+Y_{2,bq^{12}}^{-1}.
\end{split}
\end{equation*}
From condition (I), $a_1q^{-3\lambda_1+3}=q^7(a_2q^{\lambda_2-1})$. So one $Y_{1,b}^{-1}$ can only appear in $\chi_q(V_2(a_2q^{\lambda_2-1}))$, and so $(u_{1,b}(m_2')<0\Rightarrow b = a_1q^{-3\lambda_1+3}=q^7(a_2q^{\lambda_2-1}))$. As a consequence $v_{1,a_1q^{-3\lambda_1}}(m'm^{-1})\geq 1$. From the above explicit description of the $\chi_q(V_2(b))$, for all 
$$m''\in \mathcal{M}(L(m^{\rightarrow (1)})(\mathcal{M}(L(m^{\rightarrow (2)}))-\{m^{\rightarrow (2)}\}),$$ 
if $v_{1,a_1q^{-3\lambda_1}}(m''m^{-1})=0$ then $$\prod_{l\geq 0}Y_{1,a_1q^{-3\lambda_1-3+6l}}^{u_{1,a_1q^{-3\lambda_1+3+6l}}(m'')}=Y_{1,a_1q^{-3\lambda_1-3}}^{\epsilon}Y_{1,a_1q^{3-3\lambda_1}}Y_{1,a_1q^{9-3\lambda_1}}\cdots Y_{1,a_1q^{3\lambda_1-3}},$$ 
where $\epsilon \in\{0,1\}$. In particular we can prove as for property (2) of Lemma \ref{proofspe} that $v_{1,a_1q^{-3\lambda_1}}(m'm^{-1})\geq 1$ implies $v_{1,a_1q^{3\lambda_1}}(m'm^{-1})\geq 1$, contradiction.\qed

\begin{lem}\label{proofspegdeux} Let $m$ be a dominant monomial satisfying condition (II) of Theorem \ref{class}. Then $L(m)$ is special.\end{lem}

\demo It follows from Lemma \ref{geneun} that for $m'\in\mathcal{M}(L(m))$, if $v_2(m'm^{-1})=0$ then $(m')^{\rightarrow (2)}$ is of the form 
$$(Y_{2,a_2q^{\lambda_2-1}}Y_{2,a_2q^{\lambda_2-3}}\cdots Y_{2,a_2q^{1-\lambda_2}})Y_{2,a_2q^{1-\lambda_2-2}}\cdots Y_{2,a_2q^{1-\lambda_2-2R}},$$ 
where $R\geq 0$  (from condition (II) we have $a_2q^{1-\lambda_2-2}=q^5(a_1q^{3\lambda_1-3})$). So we can use the proof of property (2) of Lemma \ref{proofspe}.\qed

Lemma \ref{proofspeg} and Lemma \ref{proofspegdeux} combined with corollary \ref{simplification} imply Theorem \ref{fora} for type $G$.

\subsection{Types $C$, $D$ and $F_4$} In this subsection we prove theorem \ref{forcd}.

From corollary \ref{simplificationcd}, it suffices to consider the condition (II).

Type $C$ : as $\lambda_n = 0$ and $\Glie_{\{1,\cdots,n-1\}}$ is of type $A_{n-1}$, it follows from (1) of Lemma \ref{proofspe} that the monomials $m'\in\mathcal{M}(L(m))$ satisfying $v_n(m'm^{-1}) > 0$ are right-negative and so not dominant. For the monomials $m'\in\mathcal{M}(L(m))$ satisfying $v_n(m'm^{-1}) = 0$, we can use (2) of Lemma \ref{proofspe} and Lemma \ref{depart}.

Type $D$ : as $a_n = a_{n-1}$ and $\lambda_n = \lambda_{n-1}$, all monomials in the set
$$m^{\rightarrow (I-\{n-1,n\})}\mathcal{M}(L(m^{\rightarrow (n)}))\mathcal{M}(L(m^{\rightarrow (n - 1)}))$$ 
are right-negative. Moreover we can prove as (1) of Lemma \ref{proofspe} that for $i=n-1$ or $i=n$, $v_i(m'm^{-1}) > 0$ implies $v_{i,a_iq^{\lambda_i}}(m'm^{-1}) > 0$, and so $m'$ is right-negative. For the monomials $m'\in\mathcal{M}(L(m))$ satisfying $v_{n-1}(m'm^{-1}) = v_n(m'm^{-1})= 0$, we can use (2) of Lemma \ref{proofspe} and Lemma \ref{depart}.

Type $F_4$ : the proof is analog to type $C$ by using Lemma \ref{proofspebdeux} for $\Glie_{\{1,2,3\}}$ of type $B_3$.\qed

\section{Applications and further possible developments}\label{lastsec}

\subsection{Jacobi-Trudi determinants and Nakai-Nakanishi conjecture} In \cite[Conjecture 2.2]{nn1} Nakai-Nakanishi conjectured for classical types that the Jacobi-Trudi determinant is the $q$-character of a certain finite dimensional representation of the corresponding quantum affine algebra. This determinant can be expressed in terms of tableaux (see \cite{br} for type $A$, \cite{kos} for type $B$, and \cite{nn1, nn2, nn3} for general classical type). The cases considered in \cite{nn1} include all minimal affinizations for type $A$, and for type $B$ many minimal affinizations (but for example not the fundamental representations $V_n(a)$).

As an application of the present paper, we prove this conjecture for minimal affinizations of type $A$ and $B$ considered in \cite[Conjecture 2.2]{nn1} (see the introduction for previous results). Indeed it can be checked for type $A$ and $B$ that the tableaux expression is special and canceled by screening operators, and so is given by the Frenkel-Mukhin algorithm (see the proofs bellow; this fact was first announced and observed in some cases in \cite[Section 2.3, Rem. 1]{nn1}). So from \cite{Fre2}, Theorem \ref{fora} proved in the present paper implies that the $q$-character of a considered minimal affinization is necessarily equal to the corresponding expression. 


\begin{thm}\label{nnconj} For $\Glie$ of type $A$, $B$, the $q$-character of a minimal affinization considered in \cite[Conjecture 2.2]{nn1} is given by the corresponding Jacobi-Trudi determinant.\end{thm}

This result is coherent with the thin property proved in this paper.

With the same strategy, representations more general than minimal affinizations, and types $C$, $D$, will be discussed in a separate publication.

Let us recall the tableaux expression of the Jacobi-Trudi determinant and give the proof of theorem \ref{nnconj}. We treat the type $B$ (the proof for type $A$ is more simple).

We recall that a partition $\lambda = (\lambda_1,\lambda_2,\cdots)$ is a sequence of weakly decreasing non-negative integers with finitely many non-zero terms. The conjugate partition is denoted by $\lambda' = (\lambda_1',\lambda_2',\cdots)$. For $\lambda,\mu$ two partitions, we say that $\mu\subset \lambda$ if for all $i\geq 0$, $\lambda_i\geq \mu_i$. For $\mu\subset\lambda$, the corresponding skew diagram is 
$$\lambda/\mu = \{(i,j)\in\NN\times\NN|\mu_i + 1\leq j\leq \lambda_i\} = \{(i,j)\in\NN\times\NN|\mu_j' + 1\leq i\leq \lambda_j'\}.$$ 
We suppose in the following that $d(\lambda/\mu) \leq n$ where $d(\lambda/\mu)$ is the length of the longest column of $\lambda/\mu$, and that $\lambda/\mu$ is connected (i.e. $\mu_i + 1\leq \lambda_{i+1}$ if $\lambda_{i+1}\neq 0$).

Let $\mathbf B = \{ 1,\cdots,n,0,\overline{n},\cdots, \overline{1}\}$. 
We give the ordering $\prec$ on the set $\mathbf B$
by
\begin{equation*}
  1 \prec 2 \prec \cdots \prec n \prec 0 \prec \overline{n} \prec\cdots \prec \overline{2}\prec \overline{1}.
\end{equation*}
As it is a total ordering, we can define the corresponding maps $\text{succ}$ and $\text{prec}$.
For $a\in\CC^*$, let
{\allowdisplaybreaks
\begin{equation*}
\begin{aligned}[c]
    & \ffbox{1}_a = Y_{1,a}, 
    \\
    & \ffbox{i}_a = Y_{i-1,aq^{2i}}^{-1} Y_{i,aq^{2(i-1)}} \qquad(2 \le i \le n-1),
    \\
    & \ffbox{n}_a
     = Y_{n-1,aq^{2n}}^{-1} Y_{n,aq^{2n - 1}}Y_{n,aq^{2n - 3}},
    \\
    & \ffbox{0}_a = Y_{n,aq^{2n+1}}^{-1} Y_{n,aq^{2n-3}},
    \\
    & \ffbox{\overline{n}}_a = Y_{n-1,aq^{2n - 2}} Y_{n,aq^{2n+1}}^{-1}Y_{n,aq^{2n-1}}^{-1},
    \\
    & \ffbox{\overline{i}}_a = Y_{i-1,aq^{4n - 2i - 2}} Y_{i,aq^{4n - 2i}}^{-1}
    \qquad(2 \le i \le n-1),    \\
    & 
    \ffbox{\overline{1}}_a = Y_{1,aq^{4n - 2}}^{-1}.
\end{aligned}
\end{equation*}}
Observe that we have
$$\chi_q(V_1(a)) = \ffbox{1}_a  + \ffbox{2}_a + \cdots + \ffbox{n}_a + \ffbox{0}_a + \ffbox{\overline{n}}_a + \ffbox{\overline{n-1}}_a + \cdots + \ffbox{\overline{1}}_a.$$
For $T = (T_{i,j})_{(i,j)\in\lambda/\mu}$ a tableaux of shape $\lambda/\mu$ with coefficients in $\mathbf B$, let 
$$m_{T,a} = \prod_{(i,j)\in\lambda/\mu} \ffbox{T_{i,j}}_{aq^{4(j-i)}}\in\Yim.$$
Let $\text{Tab}(B_n,\lambda/\mu)$ be the set of tableaux of shape $\lambda/\mu$ with coefficients in $\mathbf B$ satisfying the two conditions :

$T_{i,j}\preceq T_{i,j+1}$ and $(T_{i,j},T_{i,j+1})\neq (0,0)$,

$T_{i,j}\prec T_{i+1,j}$ or $(T_{i,j},T_{i+1,j}) = (0,0)$.

The tableaux expression of the Jacobi-Trudi determinant \cite{kos, nn1} is :
$$\chi_{\lambda/\mu,a} = \sum_{T\in \text{Tab}(B_n,\lambda/\mu)} m_{T,a} \in \Yim .$$

For a monomial $m$, we denote $(m)^{\pm} = \prod_{\{i\in I, a\in\CC^*|\pm u_{i,a}(m) > 0\}}Y_{i,a}^{u_{i,a}(m)}$ the negative and the positive part of $m$.

We say that $(m)^-$ is partly canceled by $(m')^+$ if there is $i\in I$ and $a\in\CC^*$ such that $u_{i,a}((m)^-) = - u_{i,a}((m)^+)\neq 0$.

\begin{lem}\label{cancelboite} Let $T\in \text{Tab}(B_n,\lambda/\mu)$ and $a\in\CC^*$. Let $(i,j) \neq (i',j')\in \lambda/\mu$, $\alpha = T_{i,j}$ and $\beta = T_{i',j'}$. If $(\ffbox{\alpha}_{aq^{4(j-i)}})^-$ is partly canceled by $(\ffbox{\beta}_{aq^{4(j'-i')}})^+$, then $(i,j)=(i'+1,j')$ or $((i,j) = (i'+1,j'+1)\text{ and }T_{i,j} = \overline{T_{i',j'}} = \overline{n})$.\end{lem}

\demo We study different cases : 

Case (1) : $2\preceq\alpha\preceq n$ and $1\preceq\beta\preceq n-1$. We have $\alpha = \beta + 1$ and $q^{4(j-i) + 2\alpha} = q^{4(j'-i')+2(\beta - 1)}$. So $j'-i' = (j - i) + 1$. If $j < j'$, we have $i\leq i'$ and so $T_{i,j}\preceq T_{i',j'}$, contradiction. So $j\geq j'$ and $i > i'$. There is $((i_r,j_r))_{1\leq r\leq R}\in (\lambda/\mu)^R$ such that $(i_0,j_0)=(i',j')$ and $(i_R,j_R)=(i,j)$ and ($(i_{r+1},j_{r+1})=(i_r+1,j_r)$ or $(i_{r+1},j_{r+1})=(i_r,j_r+1)$). Let $T_r = T_{i_r,j_r}$. As $(i_{r+1},j_{r+1})=(i_r+1,j_r)$ implies $n\succeq T_{r+1}\succ T_r$, we have $T_R\succeq T_1 + (i-i')$, and so $(i,j)=(i'+1,j')$.

Case (2) : $\overline{n-1}\preceq \alpha\preceq \overline{1}$ and $\overline{n}\preceq \beta\preceq \overline{2}$. Analog to case (1).

Case (3) : $2\preceq\alpha\preceq n$ and $\overline{n}\preceq \beta\preceq \overline{2}$. As $$\ffbox{\alpha}_{aq^{4(j-i)}}\in\ZZ[Y_{k,aq^{2k-2+4r}}]_{k\leq n-1,r\in\ZZ}\times\ZZ[Y_{n,aq^{2r}}]_{r\in\ZZ},$$ 
and 
$$\ffbox{\beta}_{aq^{4(j'-i')}}\in\ZZ[Y_{k,aq^{2k+4r}}]_{k\leq n-1,r\in\ZZ}\times\ZZ[Y_{n,aq^{2r}}]_{r\in\ZZ},$$ 
we have a contradiction.

Case (4) : $\overline{n-1}\preceq \alpha\preceq \overline{1}$ and $1\preceq\beta\preceq n-1$. Analog to case (3).

Case (5) : $\alpha = 0$ and $\beta = n$. We have $q^{4(j-i) + 2n+1} = q^{4(j'-i')+2n - 3}$. So $j'-i' = (j - i) + 1$. As in case (1), we have $j\geq j'$. So $i > i'$. Consider $(i_r,j_r)$, $T_r$ as in case (1). If $i \geq i' + 2$, there is $r_1<r_2$ such that $i_{r_1+1} = i_{r_1} +1$ and $i_{r_2+1} = i_{r_2} + 1$. We have $T_{r_1} = T_{r_1+1} = 0$ or $T_{r_2} = T_{r_2+1} = 0$. So there is $(p,q)\in \lambda/\mu$ such that $(p,q+1), (p+1,q+1)\in\lambda/\mu$ and $T_{p,q+1}=T_{p+1,q+1}=0$ and $T_{p,q}=n$. So $(p+1,q)\in\lambda/\mu$ and $T_{p+1,q}=n$, contradiction.

Case (6) : $\alpha = 0$ and $\beta = 0$. We have $q^{4(j-i) + 2n+1} = q^{4(j'-i')+2n - 3}$ and we can conclude as in case (5).

Case (7) : $\alpha = \overline{n}$ and $\beta = 0$. We have $q^{4(j-i) + 2n+1} = q^{4(j'-i')+2n - 3}$. So $j'-i' = (j - i) + 1$. As in case (1) we have $j\geq j'$. So $i > i'$. If $j > j'$, as in case (5) we get $(p,q)\in \lambda/\mu$ such that $(p+1,q), (p+1,q+1)\in\lambda/\mu$ and $T_{p,q}=T_{p+1,q}=0$ and $T_{p+1,q+1}=\overline{n}$. So $(p,q+1)\in\lambda/\mu$ and $T_{p,q+1}=\overline{n}$, contradiction.

Case (8) : $\alpha = \overline{n}$ and $\beta = n$. We have $q^{4(j-i) + 2n+1} = q^{4(j'-i')+2n - 3}$ or $q^{4(j-i) + 2n - 1} = q^{4(j'-i')+2n - 1}$. In the first case $j'-i' = (j - i) + 1$. As above we have $j\geq j'$. So $i > i'$. Consider $(i_r,j_r), T_r$ as in case (1). If there is $r$ such that $((i_r,j_r),(i_{r+1},j_{r+1}),(i_{r+2},j_{r+2}))=((i_r,j_r),(i_r,j_r+1),(i_r+1,j_r+1))$, we have necessarily $(T_r,T_{r+1},T_{r+2}) = (n,0,\overline{n})$. So $i'=i_r$ and $i=i_r+1=i'+1$. We can treat in the same way the situation where there is $r$ such that $((i_r,j_r),(i_{r+1},j_{r+1}),(i_{r+2},j_{r+2}))=((i_r,j_r),(i_r+1,j_r),(i_r+1,j_r+1))$. In the second case $j'-i' = (j - i)$. As above we have $j > j'$ and $i=i'+1$.\qed

\begin{lem} Let $T_0 = (i - \mu_j')_{(i,j)\in\lambda/\mu}$. Then $T_0\in\text{Tab}(B_n,\lambda/\mu)$ and $m_{T_0,a}$ is the unique dominant monomial of $\chi_{\lambda/\mu,a}$.\end{lem}

\demo First it is clear that $T_0 \in \text{Tab}(B_n,\lambda/\mu)$ and that $m_{T_0,a}$ is dominant. Consider $T\in\text{Tab}(B_n,\lambda/\mu)$ such that $T_0\neq T$. So there is $(i,j)\in \lambda/\mu$ satisfying the property
\begin{equation}\label{propan}
(i = \mu_j' + 1\text{ and }T_{i,j}\neq 1)\text{ or }(i\neq \mu_j' + 1\text{ and }T_{i,j}\neq \text{succ}(T_{i-1,j})).
\end{equation} 
From lemma \ref{cancelboite} the negative part of the box corresponding to $(i,j)$ is not canceled in $m_{T,a}$ (in the case (8) of lemma \ref{cancelboite}, the negative part of the box can only be partly canceled).\qed

\begin{lem}\label{thin} For all $T\in \text{Tab}(B_n,\lambda/\mu)$, $a\in\CC^*$, the monomial $m_{T,a}$ is thin.\end{lem}

\demo Let $(i,j)\neq (i',j')\in \lambda/\mu$, $\alpha = T_{i,j}$ and $\beta = T_{i',j'}$. We suppose that $(\ffbox{\alpha}_{aq^{4(j-i)}})^+ = (\ffbox{\beta}_{aq^{4(j'-i')}})^+ \neq 0$. We study different cases (by symmetry we can suppose $\alpha\preceq \beta$) : 

Case (1) : $1\preceq\alpha\preceq \beta \preceq n-1$. We have $\alpha = \beta$ and $q^{4(j-i) + 2(\alpha-1)} = q^{4(j'-i')+2(\beta - 1)}$. So $j'-i' = (j - i)$. If $j < j'$, we have $i > i'$ and so $T_{i,j}\prec T_{i',j'}\preceq n-1$, contradiction. In the same way we get a contradiction for $j > j'$.

Case (2) : $\overline{n}\preceq \alpha\preceq \beta\preceq \overline{2}$. Analog to case (1).

Case (3) : $1\preceq\alpha\preceq n-1$ and $\overline{n}\preceq \beta\preceq \overline{2}$. Analog to case (3) of lemma \ref{cancelboite}.

Case (4) : $\alpha = n$ and $\beta = 0$. We have $q^{4(j-i) + 2n - 3} = q^{4(j'-i')+2n - 3}$. So $j'-i' = j - i$. As above, we have $j < j'$. So $i < i'$. We can conclude as in case (5) of lemma \ref{cancelboite}.

Case (5) : $\alpha =\beta = 0$. We have $q^{4(j-i) + 2n - 3} = q^{4(j'-i')+2n - 3}$. So $j'-i' = j - i$. If $j\neq j'$, we get $(p,q)\in\lambda/\mu$ such that $(p,q+1)\in\lambda/\mu$ and $T_{p,q}=T_{p,q+1}=0$, contradiction.

Case (6) : $\alpha = \beta = n$. We have $q^{4(j-i) + 2n - 3} = q^{4(j'-i')+2n - 3}$ or $q^{4(j-i) + 2n - 1} = q^{4(j'-i')+2n - 1}$. In both cases $j'-i' = j - i$ and we get a contradiction as in case (1).\qed

Finally we can conclude the proof of theorem \ref{nnconj} :

\begin{lem} We have $\chi_{\lambda/\mu,a}\in\text{Im}(\chi_q)$.\end{lem}

In the proof we will need the following partial ordering defined on $\text{Tab}(B_n,\lambda/\mu)$ : for $T,T'\in\text{Tab}(B_n,\lambda/\mu)$ we set :
$$T\preceq T'\Leftrightarrow (\forall (i,j)\in\lambda/\mu, T_{i,j}\preceq T_{i,j}').$$
Also by convention for any $\alpha\in {\bf B}$, $T_{i,j}\neq \alpha$ means ($(i,j)\in\lambda/\mu \Rightarrow T_{i,j} \neq \alpha$).

\demo Let $\alpha\in I$. We want to give a decomposition of $\chi_{\lambda/\mu}$ as in proposition \ref{jdecomp} for $J = \{\alpha\}$. From Lemma \ref{thin}, the $L_{\alpha}(M)$ that should appear in this decomposition are thin. It suffices to prove that the set $\text{Tab}(B_n,\lambda/\mu)$ is in bijection with a disjoint union of sets $\mathcal{M}(L_{\alpha}(M))$ via $T\mapsto m_{T,a}$. 

First suppose that $\alpha\leq n - 1$. Let $\mathcal{M}_{\alpha}$ be the set of tableaux $T\in\text{Tab}(B_n,\lambda/\mu)$ such that for any $(i,j)\in \lambda/\mu$ :
$$T_{i,j} = \alpha + 1 \Rightarrow ((i-1,j)\in\lambda/\mu\text{ and }T_{i-1,j} = \alpha),$$
$$T_{i,j} = \overline{\alpha}  \Rightarrow ((i-1,j)\in\lambda/\mu\text{ and }T_{i-1,j} = \overline{\alpha + 1}).$$
Then by Lemma \ref{cancelboite}, $\mathcal{M}_{\alpha}$ corresponds to all $\alpha$-dominant monomials appearing in $\chi_{\lambda/\mu,a}$. For $T\in\mathcal{M}_\alpha$, let $\tilde{T}$ be the tableaux defined in the following way. For $(i,j)\in\lambda/\mu$ :

if $T_{i,j} = \alpha\text{ and }T_{i+1,j} \neq \alpha + 1$, we set $\tilde{T}_{i,j} = \alpha + 1$,

if $T_{i,j} = \overline{\alpha + 1}\text{ and }T_{i+1,j} \neq \overline{\alpha}$, we set $\tilde{T}_{i,j} = \overline{\alpha}$,

otherwise we set $\tilde{T}_{i,j} = T_{i,j}$.

\noindent Then $\tilde{T}\in\text{Tab}(B_n,\lambda/\mu)$. For $T\in\mathcal{M}_\alpha$, we define :
$$\mathcal{M}_\alpha(T) = \{T'\in\text{Tab}(B_n,\lambda/\mu)|T\preceq T' \preceq \tilde{T}\}.$$
Then by Lemma \ref{aidesldeux} we have
$$L_\alpha(m_{T,a}) = \sum_{T'\in\mathcal{M}_\alpha(T)} m_{T',a},$$
and $(\mathcal{M}_\alpha(T))_{T\in\mathcal{M}_\alpha}$ defines a partition of $\text{Tab}(B_n,\lambda/\mu)$.

Now we treat the case $\alpha = n$.  Let $\mathcal{M}_n$ be the set of tableaux $T\in\text{Tab}(B_n,\lambda/\mu)$ such that for any $(i,j)\in \lambda/\mu$ :
$$T_{i,j} = 0 \Rightarrow ((i-1,j)\in\lambda/\mu\text{ and }T_{i-1,j} \in \{0,n\}),$$
$$T_{i,j} = \overline{n}  \Rightarrow ((i-1,j-1)\in\lambda/\mu\text{ and } T_{i-1,j-1} = n).$$
By definition of skew diagram, the last condition implies that 
$$(T_{i,j} = \overline{n} \Rightarrow ((i-1,j),(i,j-1)\in\lambda/\mu\text{ and }T_{i-1,j} \in\{0,n\}\text{ and }T_{i,j-1}\in\{0,\overline{n}\})).$$ 
This can be rewritten :
$$T_{i,j} = \overline{n} \Rightarrow \begin{pmatrix}T_{i-1,j-1} & T_{i,j-1}\\T_{i-1,j} & T_{i,j} \end{pmatrix}\in\{\begin{pmatrix} n & 0\\ n & \overline{n}\end{pmatrix},\begin{pmatrix} n & 0\\ 0 & \overline{n}\end{pmatrix},\begin{pmatrix} n & \overline{n}\\ n & \overline{n}\end{pmatrix},\begin{pmatrix} n & \overline{n}\\ 0 & \overline{n}\end{pmatrix}\}.$$
Then by Lemma \ref{cancelboite}, $\mathcal{M}_n$ corresponds to all $n$-dominant monomials appearing in $\chi_{\lambda/\mu,a}$. For $T\in\mathcal{M}_n$, let $\tilde{T}$ be the tableaux defined in the following way. For $(i,j)\in\lambda/\mu$ :

if $T_{i,j} = n$ and $T_{i+1,j+1}\neq\overline{n}$ and $T_{i+1,j}\neq 0$ and $T_{i+1,j}\neq \overline{n}$, we set $\tilde{T}_{i,j} = \overline{n}$,

if $T_{i,j} = n$ and $T_{i+1,j+1}\neq\overline{n}$ and $T_{i+1,j}\in\{0,\overline{n}\}$, we set $\tilde{T}_{i,j} = 0$,

if $T_{i,j} = 0$ and $T_{i+1,j}\neq 0$ and $T_{i+1,j} \neq \overline{n}$, we set $\tilde{T}_{i,j} = \overline{n}$.

otherwise we set $\tilde{T}_{i,j} = T_{i,j}$.

\noindent Then $\tilde{T}\in\text{Tab}(B_n,\lambda/\mu)$. For $T\in\mathcal{M}_n$, we define :
$$\mathcal{M}_\alpha(T) = \{T'\in\text{Tab}(B_n,\lambda/\mu)|T\preceq T' \preceq \tilde{T}\}.$$
Then by Lemma \ref{aidesldeux}, we have
$$L_n(m_{T,a}) = \sum_{T'\in\mathcal{M}_\alpha(T)} m_{T',a},$$
and $(\mathcal{M}_n(T))_{T\in\mathcal{M}_n}$ defines a partition of $\text{Tab}(B_n,\lambda/\mu)$.\qed


\subsection{General quantum affinizations} The quantum affinization $\U_q(\hat{\Glie})$ of a quantum Kac-Moody algebra $\U_q(\Glie)$ is defined with the same generators and relations as the Drinfeld realization of quantum affine algebras, but by using the generalized symmetrizable Cartan matrix of $\Glie$ instead of a Cartan matrix of finite type. The quantum affine algebra, quantum affinizations of usual quantum groups, are the simplest examples and have the particular property of being also quantum Kac-Moody algebras. In general these algebras are not a quantum Kac-Moody algebra. In \cite{mi, Naams, her04}, the category $\mathcal{O}$ of integrable representations is studied. For regular quantum affinizations (with a linear Dynkin diagram), one can define analogs of minimal affinizations by using properties (I) and (II) of Theorem \ref{class}. 

For example let us consider the type $B_{n,p}$ ($n\geq 2$, $p\geq 2$) corresponding to the Cartan matrix $(C_{i,j})_{1\leq i,j\leq n}$ defined as the Cartan matrix of type $B_n$ except that we replace $C_{n,n-1} = -2$ by $C_{n,n-1} = -p$. Then one can prove exactly as for lemma \ref{proofspebdeux} that (an analog of Theorem \ref{racourc} is proved by using \cite[Lemma 5.10]{her04}):

\begin{thm}\label{moregen} Let $\Glie$ be of type $B_{n,p}$. Then if $m$ satisfies property (I) (resp. (II)), then $L(m)$ is antispecial (resp. special).\end{thm}

So the analog of the Frenkel-Mukhin algorithm works for these modules and as an application it should be possible to get additional results for this class of special modules (see also section \ref{classgen} bellow).

\subsection{Multiparameter $T$-systems} 

The special property of Kirillov-Reshetikhin modules allows to prove a system of induction relations involving $q$-characters of Kirillov-Reshetikhin modules called $T$-system (see \cite{Nad} for the simply-laced cases and \cite{her06} for the general case). Indeed for $i\in I$, $k\geq 1$, $a\in\CC^*$ define the $\U_q(\Lo\Glie)$-module : 
$$S_{r,a}^{(i)}=({\bigotimes}_{\{(j,k)|C_{j,i}<0\text{ , }1\leq k\leq -C_{i,j}\}}W_{-C_{j,i}+E(r_i(r-k)/r_j),aq_j^{-(2k-1)/C_{i,j}}}^{(j)}).$$

\begin{thm}[The $T$-system]\label{tsyst} Let $a\in\CC^*, k\geq 1, i\in I$. Then we have :
$$\chi_q(W_{k,a}^{(i)})\chi_q(W_{k,aq_i^2}^{(i)})=\chi_q(W_{k+1,a}^{(i)})\chi_q(W_{k-1,aq_i^2}^{(i)})
+\chi_q(S_{k,a}^{(i)}).$$ 
\end{thm}

By analogy, the results of the present paper (special property of minimal affinizations of type $A,B,G$) should lead to systems of induction relations involving $q$-characters of minimal affinizations (multiparameter $T$-systems). Let us look at an example. Let $\Glie = sl_3$. Then we have the following relation :
$$\chi_q(L(X_{3,q^2}^{(1)}X_{2,q^8}^{(2)}))\chi_q(L(X_{3,q^4}^{(1)}X_{2,q^{10}}^{(2)}))$$
$$= \chi_q(L(X_{4,q^3}^{(1)}X_{2,q^{10}}^{(2)}))\chi_q(L(X_{2,q^3}^{(1)}X_{2,q^8}^{(2)}))+\chi_q(L(X_{6,q^6}^{(2)}))\chi_q(L(X_{1,q^9}^{(2)})).$$
Let us give the idea of the proof for this example : as a $q$-character is characterized by the multiplicity of his dominant monomials \cite{Fre2}, it suffices to compare dominant monomial of both side. By using the process described in remark \ref{process}, Theorem \ref{racourc} and arguments of \cite{her06}, we get the following results :

\noindent The dominant monomials of $\chi_q(L(X_{3,q^2}^{(1)}X_{2,q^8}^{(2)})\otimes L(X_{3,q^4}^{(1)}X_{2,q^{10}}^{(2)}))$ are :

$1_0 1_2^2 1_4 2_5 2_7 2_9^2 2_11$,
$1_0 1_2 2_3 2_5 2_7 2_9^2 2_11$,
$2_1 2_3 2_5 2_7 2_9^2 2_11$,
$1_0 1_2 1_4^2 1_6 1_{10} 2_7 2_9$,
\\$1_0 1_2^2 1_4 1_{10} 2_5 2_7 2_9$,
$1_0 1_2 1_{10} 2_3 2_5 2_7 2_9$,
$1_{10} 2_1 2_3 2_5 2_7 2_9$.

\noindent The dominant monomials of $\chi_q(L(X_{4,q^3}^{(1)}X_{2,q^{10}}^{(2)})\otimes L(X_{2,q^3}^{(1)}X_{2,q^8}^{(2)}))$ are :

$1_0 1_2^2 1_4 2_5 2_7 2_9^2 2_{11}$,
$1_0 1_2 2_3 2_5 2_7 2_9^2 2_{11}$,
$1_0 1_2^2 1_4^2 1_6 1_{10} 2_7 2_9$,
$1_0 1_2^2 1_4 1_{10} 2_5 2_7 2_9$,
\\$1_0 1_2 1_{10} 2_3 2_5 2_7 2_9$.

\noindent The dominant monomials of $\chi_q(L(X_{3,q^2}^{(1)}X_{2,q^8}^{(2)})\otimes L(X_{3,q^4}^{(1)}X_{2,q^{10}}^{(2)}))$ are :

$2_1 2_3 2_5 2_7 2_9^2 2_{11}$,
$1_{10} 2,q 2_3 2_5 2_7 2_9$.

\noindent We can conclude as the multiplicity of all these monomials is $1$.

\subsection{Alternative method for the classification of minimal affinizations}\label{classgen}

We explain how to prove certain classification results (included in Theorem \ref{class}). The proofs here are written in the context of the paper and could be a general uniform strategy for other quantum affinizations. Moreover we get some new refined results on the involved $q$-characters.

\begin{prop}\label{qseg} Let $L(m)$ be a minimal affinization of $V(\lambda)$. Then for all $i\in I$, there is $a_i\in\CC^*$ such that $m^{\rightarrow (i)}=X_{a_i,\lambda_i}^{(i)}$.\end{prop}

\demo For $\lambda_i\leq 1$ it is clear. Suppose that $\lambda_i\geq 2$ and that $m^{\rightarrow (i)}$ in not of this form. Note that $\lambda - \alpha_i\in P^+$. It follows from Lemma \ref{depart} with $J=\{i\}$ and Proposition \ref{aidesldeux} that 
$$\text{dim}((L(m))_{\lambda -\alpha_i})=\text{dim}((L_i(m^{\rightarrow (i)}))_{(\lambda_i-2)\Lambda_i}\geq 2.$$ 
Let $a\in\CC^*$ and $M=m((m)^{\rightarrow (i)})^{-1}X_{\lambda_i,a}^{(i)}$. $L(M)$ is an affinization of $V(\lambda)$. It follows from Lemma \ref{depart} with $J=\{i\}$ that $\text{dim}((L(M))_{\lambda-\alpha_i})=1$ so $m_{\lambda - \alpha_i}(L(M)) < m_{\lambda - \alpha_i}(L(m))$. Moreover as $(m)^{\rightarrow (I-\{i\})}=(M)^{\rightarrow (I-\{i\})}$, it follows from Lemma \ref{depart} with $J=I-\{i\}$ that for $\mu\in \lambda - {\sum}_{j\neq i}\NN \alpha_j$ we have $\text{dim}((L(M))_{\mu})=\text{dim}((L_{I-\{i\}}(m^{\rightarrow (I-\{i\})}))_{\mu})=\text{dim}((L(m))_{\mu})$ and so $m_{\mu}(L(M))=m_{\mu}(L(m))$. As $\mu\leq \lambda$ implies $\mu=\lambda$ or $\mu\leq \lambda -\alpha_i$ or $\mu\in \lambda - {\sum}_{j\neq i} \NN \alpha_j$, we have $[L(M)]<[L(m)]$, contradiction.\qed

In the following for $L(m)$ a minimal affinization and for $i\in I$ such that $\lambda_i\neq 0$, $a_i\in\CC^*$ denotes the complex number introduced in this Proposition \ref{qseg}.

Let $\Glie=sl_{n+1}$ ($n\geq 2$) and $\lambda=\lambda_1\Lambda_1+ \lambda_n\Lambda_n$ ($\lambda_1,\lambda_n\geq 1$). For $\mu=\alpha_1+\alpha_2+\cdots +\alpha_n$, we have $\text{dim}((V(\lambda))_{\lambda-\mu}) = n$. Let $m=X_{\lambda_1,a_1}^{(1)}X_{\lambda_n,a_n}^{(n)}$. If $L(m)$ is a minimal affinization of $V(\lambda)$ then $\text{dim}((L(m))_{\lambda-\mu}) = n$. For $0\leq h\leq n$ denote 
$$m_h = m \prod_{i=1,\cdots,h}A_{i,a_1q^{\lambda_1+i-1}}^{-1}\prod_{i=h+1,\cdots,n}A_{i,a_nq^{\lambda_n+n-i}}^{-1}.$$ 
We have different cases :

\begin{enumerate}
\item $a_1/a_n\notin\{q^{\pm(\lambda_1+\lambda_n + n-1)},q^{\lambda_n-\lambda_1+n-1},q^{\lambda_n-\lambda_1+n-3},\cdots,q^{\lambda_n-\lambda_1-n+1}\}$. \\From remark \ref{process}, the $n+1$ monomials $m_h$ for $0\leq h\leq n$ appear in $\chi_q(L(m))$ and are distinct. So $\text{dim}((L(m))_{\lambda-\mu})\geq n+1$ and $L(m)$ is not a minimal affinization of $V(\lambda)$.

\item $a_1/a_n = q^{\lambda_n- \lambda_1+n+1-2H}$ with $1\neq H\leq n$. Then $m_H = m_{H-1}$. \\From remark \ref{process}, the $n-1$ distinct monomials $m_h$ for $h\notin\{H-1,H\}$ appear in $\chi_q(L(m))$ with multiplicity $1$ and $m_H$ appears in $\chi_q(L(m))$ with multiplicity $2$. So $\text{dim}((L(m))_{\lambda-\mu})\geq n + 1$ and $L(m)$ is not a minimal affinization of $V(\lambda)$.

\item $a_1/a_n = q^{\lambda_1+\lambda_n + n - 1}$. 

\item $a_n/a_1=q^{\lambda_1+\lambda_n+ n - 1}$.
\end{enumerate}

From Proposition \ref{vraiinv}, the character is the same in cases (3) and (4). So necessarily these two cases give a minimal affinization with $\chi(L(m))=\chi(V(\lambda))$. So for $\lambda_1,\lambda_n >0$ , $L(m)$ is a minimal affinization of $V(\lambda_1\Lambda_1+\lambda_n\Lambda_n)$ if and only if $m=X_{\lambda_1,a_1}^{(1)}X_{\lambda_n,a_n}^{(n)}$ with $a_1/a_n=q^{\lambda_1+\lambda_n+n-1}$ or $a_n/a_1=q^{\lambda_1+\lambda_n+n-1}$.

Now we suppose that $\Glie$ is general and consider $J\subset I$ such that $\Glie_J$ is of type $A_r$, $2\leq r\leq n$. Denote by $i,j\in J$ the two extremes nodes of $J$. We suppose that we can decompose $I=I_i\sqcup J \sqcup I_j$ such that $I_i\cup\{i\}$ and $I_j\cup\{j\}$ are connected, and $\forall k\in I_i, k'\in J-\{i\}, C_{k,k'}=0$ and $\forall k\in I_j,k'\in J-\{j\}, C_{k,k'}=0$. Observe that $I_i$ or $I_j$ may be empty and if $J$ is of type $A_2$ there is always such a decomposition. 

\begin{prop}\label{deuxcon} Let $L(m)$ be a minimal affinization of $V(\lambda)$ such that $\lambda_i,\lambda_j\geq 1$ and for $k\in J-\{i,j\}$, $\lambda_k = 0$. Then one of the two following condition holds :
\begin{equation*}
\frac{a_i}{a_j}=q_i^{\lambda_i+\lambda_j+r-1}\text{ or }\frac{a_j}{a_i}=q_i^{\lambda_i+\lambda_j+r-1}.
\end{equation*}
\end{prop}

\demo We can suppose in the proof that $q_i = q_j = q$. Suppose that $a_i/a_j\neq q^{\pm (\lambda_i+\lambda_j+r+1)}$. Note that $\lambda-\sum_{k\in J}\alpha_k\in P^+$. It follows from Lemma \ref{depart} with $J$ and the above discussion that $\text{dim}((L(m))_{\lambda-\sum_{k\in J}\alpha_k}) \geq r+1$. Let us define 
$$M = m^{\rightarrow (I_i \cup \{i\})}\tau_{q^{\lambda_i+\lambda_j+m-1}a_ia_j^{-1}}(m^{\rightarrow (\{j\}\cup I_j)}).$$ 
$L(M)$ is an affinization of $V(\lambda)$. Let us prove that $[L(M)]<[L(m)]$ (which is a contradiction). Let $\omega\leq\lambda$. If $\omega\leq \lambda - \sum_{k\in J}\alpha_k$ it follows from Lemma \ref{depart} with $J$ that $\text{dim}((L(M))_{\lambda-\sum_{k\in J}\alpha_k}) < \text{dim}((L(m))_{\lambda-\sum_{k\in J}\alpha_k})$. As for $J'\subset J$, $\lambda-\sum_{k\in J'}\alpha_k\notin P^+$ except for $J'=J$ or $J=\emptyset$, we get $m_{\lambda-\sum_{k\in J}\alpha_k}(L(M)) < m_{\lambda-\sum_{k\in J}\alpha_k}(L(m))$. Otherwise it follows from Lemma \ref{depart} that $\text{dim}((L(M))_{\mu}) = \text{dim}((L(m))_{\mu})$ as $(M)^{\rightarrow ((I_i\cup J)-\{j\})}=(m)^{\rightarrow ((I_i\cup J)-\{j\})}$ and $(M)^{\rightarrow ((I_j\cup J)-\{i\})}=\tau_{q^{\lambda_i+\lambda_j+m-1}a_ia_j^{-1}}(m^{\rightarrow ((I_j\cup J)-\{i\})})$. So $m_{\mu}(L(M)) = m_{\mu}(L(m))$.\qed

Let $\Glie$ of type $B_n$ ($n\geq 2$), $\lambda=\lambda_1\omega_1+ \lambda_n\omega_n$ ($\lambda_1,\lambda_n\geq 1$) and $\mu=\alpha_1+\alpha_2+\cdots +\alpha_n$. Let $m=X_{\lambda_1,a_1}^{(1)}X_{\lambda_n,a_n}^{(n)}$. For $0\leq h\leq n$ denote 
$$m_h = m \prod_{i=1,\cdots,h}A_{i,a_1q^{\lambda_1+i-1}}^{-1}\prod_{i=h+1,\cdots,n}A_{i,a_nq^{2\lambda_n+1+n-i}}^{-1}.$$ We have $(L(m))_{\lambda-\mu} = {\bigoplus}_{0\leq h\leq n}(L(m))_{m_h}$. Let us study the different cases :
\begin{enumerate}
\item $a_1/a_n\notin\{q^{\pm(\lambda_1+2\lambda_n + n)},q^{2\lambda_n-\lambda_1+n},q^{2\lambda_n-\lambda_1+n-2},\cdots,q^{2\lambda_n-\lambda_1-n+2}\}$. From remark \ref{process} the $n+1$ monomials $m_h$ for $0\leq h\leq n$ appear in $\chi_q(L(m))$ and are distinct. So $\text{dim}((L(m))_{\lambda-\mu})\geq n+1$.

\item $a_1/a_n = q^{2\lambda_n- \lambda_1+n+2-2H}$ with $1\neq H\leq n$. Then $m_H = m_{H-1}$. From remark \ref{process}, the $n-1$ distinct monomials $m_h$ for $h\notin\{H-1,H\}$ appear in $\chi_q(L(m))$ with multiplicity $1$ and $m_H$ appears in $\chi_q(L(m))$ with multiplicity $2$. So $\text{dim}((L(m))_{\lambda-\mu})\geq n + 1$.

\item $a_1/a_n = q^{\lambda_1+2\lambda_n + n}$. Then $\text{dim}((L(m))_{\lambda -\mu}) = n$. Indeed We see as for the proof of the point (3) of Lemma \ref{proofspe} that for $m'\in\mathcal{M}(L(m))$, if $v_1(m'm^{-1})\geq 1$ then $v_{1,a_1q^{\lambda_1}}(m'm^{-1})\geq 1$. So $m_0\notin\mathcal{M}(L(m))$ and from remark \ref{process} $m_1,\cdots,m_H$ appear in $\chi_q(L(m))$ with multiplicity $1$.

\item $a_n/a_1=q^{\lambda_1+2\lambda_n+ n}$. As in the case (3), $\text{dim}((L(m))_{\lambda -\mu}) = n$.
\end{enumerate}

From Proposition \ref{vraiinv}, the character is the same in cases (3) and (4).

\begin{prop} For $\Glie$ of type $B_n$ with $n\geq 2$ and $\lambda_1,\lambda_n >0$ , $L(m)$ is a minimal affinization of $V(\lambda_1\Lambda_1+\lambda_n\Lambda_n)$ if and only if $m=X_{\lambda_1,a_1}^{(1)}X_{\lambda_n,a_n}^{(n)}$ with $a_1/a_n=q^{\lambda_1+2\lambda_n+n}$ or $a_n/a_1=q^{\lambda_1+2\lambda_n+n}$.\end{prop}

\demo If $m'$ satisfies (1) or (2) and $m$ satisfies (3) or (4), then $\text{dim}((L(m))_{\lambda-\mu})<\text{dim}((L(m'))_{\lambda-\mu})$ and for $\lambda'\leq \lambda$ if there is $j\in I$ such that $v_j(\lambda'-\lambda)=0$ then 
\begin{equation*}
\begin{split}
\text{dim}((L(m))_{\lambda'})
=&\text{dim}((L(m'))_{\lambda'})
\\=&\text{dim}(W_{\lambda_1,1}^{(1)})_{\lambda_1\Lambda_1-\sum_{k<j}v_k(\lambda'-\lambda)\alpha_k})
\\&\times\text{dim}(W_{\lambda_n,1}^{(n)})_{\lambda_n\Lambda_n-\sum_{k>j}v_k(\lambda'-\lambda)\alpha_k}).
\end{split}
\end{equation*}
As we have the same character in situations (3) and (4), they correspond necessarily to minimal affinizations.\qed

Now we suppose that $\Glie$ is general and consider $J\subset I$ such that $\Glie_J$ is of type $B_r$, $2\leq r\leq n$. Denote by $i,j\in J$ the two extremes nodes of $J$. We suppose that we can decompose $I=I_i\sqcup J \sqcup I_j$ such that $I_i\cup\{i\}$ and $I_j\cup\{j\}$ are connected, and $\forall k\in I_i, k'\in J-\{i\}, C_{k,k'}=0$ and $\forall k\in I_j,k'\in J-\{j\}, C_{k,k'}=0$. Observe that $I_i$ or $I_j$ may be empty and if $J$ is of type $B_2$ there is always such a decomposition.

\begin{prop} Let $L(m)$ be a minimal affinization of $V(\lambda)$ such that $\lambda_i,\lambda_j\geq 1$ and for $k\in J-\{i,j\}$, $\lambda_k = 0$. Then one of the two following condition holds :
\begin{equation*}
\frac{a_i}{a_j}=q_i^{\lambda_i+2\lambda_j+r}\text{ or }\frac{a_j}{a_i}=q_i^{\lambda_i+2\lambda_j+r}.
\end{equation*}
\end{prop}

The proof is analog to the proof of the Proposition \ref{deuxcon}.

\end{document}